\documentstyle[12pt,amssymb]{amsart}
\newcommand{\eps}{\epsilon}
\newcommand{\Td}{\operatorname{Td}}
\newcommand{\NS}{\operatorname{NS}}
\newcommand{\Div}{\operatorname{Div}}

\newcommand{\herm}{\operatorname{herm}}
\newcommand{\SU}{\operatorname{SU}}
\newcommand{\Isom}{\operatorname{Isom}}
\newcommand{\U}{\operatorname{U}}
\newcommand{\und}{\underline}

\newcommand{\W}{\operatorname{W}}
\newcommand{\WH}{\operatorname{WH}}

\newcommand{\Mat}{\operatorname{Mat}}

\newcommand{\sym}{\operatorname{sym}}
\newcommand{\ch}{\operatorname{ch}}
\newcommand{\CH}{\operatorname{CH}}
\newcommand{\Cor}{\operatorname{Cor}}

\newcommand{\Sch}{\operatorname{Schr}}

\newcommand{\SL}{\operatorname{SL}}
\renewcommand{\Sp}{\operatorname{Sp}}

\newcommand{\La}{\Lambda}

\newcommand{\tors}{\operatorname{tors}}

\newcommand{\HH}{{\cal H}}
\newcommand{\KK}{{\cal K}}
\newcommand{\LL}{{\cal L}}
\newcommand{\BB}{{\cal B}}

\newcommand{\GG}{{\cal G}}

\newcommand{\CC}{{\cal C}}
\newcommand{\Ga}{\Gamma}
\newcommand{\ga}{\gamma}
\newcommand{\si}{\sigma}

\newcommand{\pg}{{\goth p}}

\newcommand{\Aut}{\operatorname{Aut}}
\newcommand{\End}{\operatorname{End}}
\newcommand{\Hom}{\operatorname{Hom}}

\newcommand{\supp}{\operatorname{supp}}

\newcommand{\ed}{\qed\vspace{3mm}}

\newcommand{\Spec}{\operatorname{Spec}}

\newtheorem{thm}{Theorem}[section]
\newtheorem{prop}[thm]{Proposition}
\newtheorem{lem}[thm]{Lemma}
\newtheorem{cor}[thm]{Corollary}
\newenvironment{ex}{\vspace{3mm} \noindent {\bf
Example.}}{\vspace{3mm}}
\theoremstyle{definition}
\newtheorem{defi}[thm]{Definition}
\newenvironment{rem}{\vspace{3mm}
\noindent {\bf Remark.}}{\vspace{3mm}}
\newenvironment{rems}{\vspace{3mm}
\noindent {\bf Remarks.}}{\vspace{3mm}}

\numberwithin{equation}{section}
      
\newcommand{\Pf}{\noindent {\it Proof}}
\newcommand{\Pic}{\operatorname{Pic}}

\renewcommand{\a}{\alpha}
\renewcommand{\b}{\beta}
\renewcommand{\S}{{\bold S}}

\newcommand{\wt}{\widetilde}
\renewcommand{\mod}{\operatorname{mod}}

\renewcommand{\AA}{{\cal A}}

\newcommand{\FF}{{\cal F}}
\renewcommand{\GG}{{\cal G}}
\newcommand{\coker}{\operatorname{coker}}

\newcommand{\Q}{{\Bbb Q}}

\newcommand{\C}{{\Bbb C}}
\newcommand{\Z}{{\Bbb Z}}

\newcommand{\la}{\lambda}
\newcommand{\PP}{\cal P}
\newcommand{\e}{\varepsilon}
\newcommand{\Br}{\operatorname{Br}}
\newcommand{\GL}{\operatorname{GL}}

\newcommand{\ov}{\overline}

\newcommand{\om}{\omega}
\newcommand{\ra}{\rightarrow}
\newcommand{\hra}{\hookrightarrow}

\newcommand{\id}{\operatorname{id}}
\newcommand{\G}{{\Bbb G}}

\newcommand{\ot}{\otimes}

\renewcommand{\O}{{\cal O}}

\newcommand{\R}{{\Bbb R}}

\newcommand{\de}{\delta}

\newcommand{\D}{{\cal D}}
\newcommand{\DD}{{\cal D}}
\newcommand{\De}{\Delta}

\newcommand{\PPic}{{\cal P}ic}
\newcommand{\sub}{\subset}
      
\title{Analogue of Weil representation for abelian schemes}
\author{A. Polishchuk}
      
\begin{document}
\maketitle
                                           
This paper is devoted to the construction of
projective actions of certain arithmetic groups
on the derived categories of coherent sheaves
on abelian schemes over a normal base $S$.
These actions are constructed by mimicking the construction
of Weil in \cite{Weil} of a projective representation
of the symplectic group $\Sp(V^*\oplus V)$ on the space of
smooth functions on the lagrangian subspace $V$.
Namely, we replace the vector space $V$ by an abelian
scheme $A/S$, the dual vector space $V^*$ by the dual abelian
scheme $\hat{A}$, and the space of functions on $V$ by
the (bounded) derived category of coherent sheaves on $A$
which we denote by $\D^b(A)$. The role of the standard
symplectic form on $V\oplus V^*$ is played by the
line bundle $\BB=p_{14}^*\PP\ot p_{23}^*\PP^{-1}$
on $(\hat{A}\times_S A)^2$ where $\PP$ is the normalized
Poincar\'e bundle on $\hat{A}\times A$. Thus, the symplectic group
$\Sp(V^*\oplus V)$ is replaced by the group
of automorphisms of $\hat{A}\times_S A$ preserving  $\BB$.
We denote the latter group by $\SL_2(A)$ (in \cite{Orlov} the
same group is denoted by $\Sp(\hat{A}\times_S A)$).
We construct an action of a central
extension of certain "congruenz-subgroup"
$\Ga(A,2)\sub\SL_2(A)$
on $\D^b(A)$.  More precisely,
if we write an element of $\SL_2(A)$
in the block form
$$\left(   \matrix   a_{11}   &  a_{12}\\  a_{21}  &  a_{22}
\endmatrix \right)$$
then the subgroup $\Ga(A,2)$
is distinguished by the condition that elements
$a_{12}\in\Hom(A,\hat{A})$ and $a_{21}\in\Hom(\hat{A},A)$
are divisible by $2$. We
construct autoequivalences $F_g$ of $\D^b(A)$ corresponding to
elements $g\in\Ga(A,2)$, such the composition
$F_{g}\circ F_{g'}$ differs from $F_{gg'}$ by tensoring
with a line bundle on $S$ and a shift in the derived category.
Thus, we get an embedding of the central extension of
$\Ga(A,2)$ by $\Z\times\Pic(S)$ into the group of
autoequivalences of $\D^b(A)$. The 2-cocycle of this central extension
is described by
structures  similar  to  the  Maslov  index  of  a triple of
lagrangian subspaces in a symplectic vector space.
However, the situation here is more complicated since the
construction of the functor $F_g$ requires a choice of a
{\it Schr\"odinger representation} of certain
Heisenberg group $G_g$ associated with $g$. The latter is a central
extension of a finite group scheme $K_g$ over $S$ by $\G_m$
such that the commutation pairing $K_g\times K_g\ra\G_m$
is non-degenerate. If the order of $K_g$ is $d^2$ then a
 Schr\"odinger   representation   of   $G_g$   is   a
representation of $G_g$ in a vector bundle of rank $d$  over
$S$ such that $\G_m\sub G_g$ acts naturally.
Any two such representations differ by tensoring with
a  line  bundle  on  $S$.  The  classical  example
due  to D.~Mumford
arises in the situation when there is a relatively ample line
bundle $L$ on an abelian scheme $\pi:A\ra S$. Then the vector
bundle $\pi_*L$ on $S$ is a Schr\"odinger representation
of the Mumford group $G(L)$ attached to $L$, this group is
a central extension of the finite group scheme $K(L)$ by
$\G_m$  where  $K(L)$  is  the  kernel  of   the   symmetric
homomorphism $\phi_L:A\ra\hat{A}$ associated with $L$.
Our  Heisenberg  groups  $G_g$ are subgroups in some Mumford
groups and there is no canonical choice of a
Schr\"odinger representation for them in general
(this ambiguity is responsible for
the   $\Pic(S)$-part   of   our   central    extension    of
$\Ga(A,2)$). Moreover,
unless $S$ is a spectrum of an algebraically closed field,
it is not even obvious
that a Schr\"odinger representation for $G_g$ exists.
Our main technical result that deals with this problem is  the
existence of a "Schr\"odinger representation" for
finite symmetric Heisenberg group schemes of odd order
established in section \ref{Schr}. Further we observe that the obstacle
to existence of such a representation is an element $\de(G_g)$ in the
Brauer group of $S$, and that the map $g\mapsto \de(G_g)$
is a homomorphism. This allows us to use the theory of arithmetic groups
to prove the vanishing of $\de(G_g)$.
                                      
When  an abelian scheme $A$ is equipped with some additional
structure (such as a symmetric line bundle) one can sometimes
extend the above action of a central extension of
$\Ga(A,2)$ on $\D^b(A)$  to  an  action  of  a
bigger group. The following two cases seem to be particularly
interesting.
Firstly, assume that a pair of line bundles on $A$ and $\hat{A}$
is given such that the composition of the corresponding
isogenies  between  $A$  and  $\hat{A}$  is  the morphism of
multiplication by some integer $N>0$. Then we can construct an
action   of   a   central   extension   of   the   principal
congruenz-subgroup $\Ga_0(N)\sub\SL_2(\Z)$ on $\D^b(A)$
(note that in this situation there is a natural embedding of $\Ga_0(N)$
into $\SL_2(A)$ but the image is not necessarily
contained in $\Ga(A,2)$). Secondly, assume that
we have a symmetric line bundle $L$ on $A$ giving rise
to  a  principal  polarization.  Then  there  is  a  natural
embedding of $\Sp_{2n}(\Z)$ into $\SL_2(A^n)$,
where $A^n$ denotes the $n$-fold fibered product over $S$,
which is an
isomorphism when $\End(A)=\Z$. In this case we construct an action of a
central extension of $\Sp_{2n}(\Z)$ on $\D^b(A^n)$.
The main point in both cases is to show the existence of
relevant Schr\"odinger representations. Both these
situations admit natural generalizations to abelian schemes
with real multiplication treated in the same way. For example,
we consider the case of an abelian scheme $A$ with
real multiplication by a ring of integers $R$ in a totally
real number field, equipped with a symmetric line bundle $L$
such that $\phi_L:A\wt{\ra}\hat{A}$ is an $R$-linear principal
polarization. Then there is an
action  of  a  central   extension   of   $\Sp_{2n}(R)$   by
$\Z\times\Pic(S)$ on $\D^b(A^n)$.
When $R=\Z$ we determine this central extension
explicitly using a presentation of
$\Sp_{2n}(\Z)$ by generators and relations. It turns out
that the $\Pic(S)$-part of this extension is induced by a
non-trivial central extension of $\Sp_{2n}(\Z)$ by $\Z/2\Z$ via
the embedding $\Z/2\Z\hra\Pic(S)$ given by an element
$(\pi_*L)^{\ot 4}\ot \ov{\om}_A^{\ot 2}$ where $\ov{\om}_A$
is the restriction of the relative canonical bundle of $A/S$
to the zero section. Also we show that the restriction of
this central extension to certain congruenz-subgroup
$\Ga_{1,2}\sub\Sp_{2n}(\Z)$ splits.
    
In the case when $S$ is
the  spectrum  of  the  algebraically   closed   field   the
constructions of this paper were developed in \cite{Weilrep} and
\cite{Orlov}. In the latter paper the $\Z$-part of the
above central extension is described in the analytic
situation. Also in  \cite{Weilrep} we have shown that
the above action of an arithmetic group on $\D^b(A)$ can
be used to construct an action of the corresponding algebraic
group over $\Q$ on the (ungraded) Chow motive of $A$.
In the present paper we extend this to the case of abelian schemes
and their relative Chow motives.
                           
Under the conjectural equivalence of $\D^b(A^n)$
with the Fukaya category of the mirror dual symplectic torus 
(see \cite{Kon}) the above projective
action of $\Sp_{2n}(\Z)$ should correspond to a natural geometric action
on the Fukaya category. The central extension by $\Z$ appears
in the latter situation due to the fact that objects are lagrangian
subvarieties together with some additional data which form a
$\Z$-torsor (see \cite{Kon}).

The paper is organized as follows. In the first two sections we study
finite   Heisenberg   group  schemes ({\it  non-degenerate
theta groups} in the terminology of \cite{MB}) and their representations.
In particular, we establish a key result (Theorem \ref{odd})
on the existence of a Schr\"odinger representation
for a symmetric Heisenberg group scheme of odd order.
In section 3 we  consider  another  analogue  of  Heisenberg
group: the central extension $H(\hat{A}\times A)$ of
$\hat{A}\times_S A$ by the Picard groupoid of line bundles on $S$.
We develop an analogue of the classical representation theory of
real Heisenberg groups  for   $H(\hat{A}\times   A)$.
Schr\"odinger representations for finite
Heisenberg groups enter into this theory as a key
ingredient for the construction of intertwining operators.
In  section  4  we  construct   a   projective   action   of
$\Ga(A,2)$ on $\D^b(A)$, in section 5 --- the corresponding
action  of an algebraic group over $\Q$ on the relative Chow
motive of $A$. In section 6 we study the group $\wt{\SL}_2(A)$,
the extension of $\SL_2(A)$ which acts on the Heisenberg groupoid.
In section 7 we
extend the action of $\Ga(A,2)$ to that of a bigger group in the situation
of abelian schemes with real multiplication. In section 8 we
study the corresponding central extension of $\Sp_{2n}(\Z)$.

All the schemes in this paper are assumed to be noetherian.
The base scheme $S$ is always assumed to be connected.
We denote
by  $\D^b(X)$  the  bounded  derived  category  of  coherent
sheaves  on  a scheme $X$. For a morphism of schemes $f:X\ra
Y$   of   finite   cohomological   dimension  we  denote  by
$f_*:\D^b(X)\ra\D^b(Y)$ (resp.  $f^*:\D^b(Y)\ra\D^b(X)$)  the
derived functor of direct (resp. inverse) image.
For any abelian scheme $A$ over $S$ we denote by $e:S\ra A$
the zero section. For an abelian scheme
$A$ (resp. morphism of abelian schemes $f$) we denote
by $\hat{A}$ (resp. $\hat{f}$) the dual abelian scheme
(resp. dual morphism).
For  every  line  bundle  $L$  on  $A$ we denote by
$\phi_L:A\ra\hat{A}$ the corresponding morphism of abelian
schemes (see \cite{Mum}).
When this is reasonable a line bundle on an abelian
scheme is tacitly assumed to be rigidified
along the zero section (one exception is provided by line
bundles pulled back from the base).
For every integer $n$ and a commutative group scheme $G$
we denote by $[n]=[n]_G:G\ra G$
the multiplication by $n$ on $G$, and  by  $G_n\sub  G$  its
kernel.
We use freely the notational  analogy  between  sheaves  and
functions writing in particular
$\FF_x=\int_Y \GG_{y,x} dy$, where $x\in X$, $y\in Y$,
$\FF\in\D^b(X)$, $\GG\in\D^b(Y\times X)$,
instead of $\FF=p_{2*}(\GG)$.
  
\section{Heisenberg group schemes}
       
Let $K$ be a finite flat group scheme over a base scheme $S$.
A   {\it  finite  Heisenberg  group  scheme}  is  a  central
extension of group schemes
\begin{equation}\label{ext}
0\ra\G_m\ra G\stackrel{p}{\ra} K\ra 0
\end{equation}
such  that  the  corresponding  commutator  form  $e:K\times
K\ra\G_m$ is a perfect pairing.
Let $A$  be an abelian scheme over $S$, $L$ be a line bundle
on $A$ trivialized along the zero section.  Then  the  group
scheme $K(L)=\{x\in A\ |\ t_x^*L\simeq L\}$ has a canonical
central extension $G(L)$ by $\G_m$ (see \cite{Mum}).
When $K(L)$ is finite, $G(L)$  is  a finite Heisenberg
group scheme.
      
A {\it symmetric} Heisenberg group scheme is an
extension $0\ra\G_m\ra G\ra K\ra0$ as above together
with an isomorphism of central extensions $G\wt{\ra} [-1]^*G$
(identical  on  $\G_m$), where $[-1]^*G$ is the pull-back of
$G$ with respect to the inversion morphism $[-1]:K\ra K$.
For  example,  if  $L$  is  a  symmetric  line  bundle on an
abelian  scheme  $A$  (i.e.  $[-1]^*L\simeq  L$)   with   a
symmetric  trivialization  along the zero section
then  $G(L)$  is  a  symmetric  Heisenberg group scheme.
      
For any integer $n$ we denote by
$G^n$ the push-forward of $G$ with respect to the morphism
$[n]:\G_m\ra\G_m$.   For   any   pair  of  central  extensions
$(G_1,G_2)$ of the same group $K$ we denote by $G_1\ot G_2$
their  sum  (given  by  the   sum   of   the   corresponding
$\G_m$-torsors). Thus, $G^n\simeq G^{\ot n}$.
Note  that  we  have  a  canonical  isomorphism  of  central
extensions
\begin{equation}\label{inver}
G^{-1}\simeq [-1]^*G^{op}
\end{equation}
where $[-1]^*G^{op}$ is the
pull-back of the opposite group  to  $G$  by  the  inversion
morphism  $[-1]:K\ra  K$. In particular, a symmetric extension
$G$ is commutative  if  and  only  if  $G^2$  is trivial.

\begin{lem}\label{mult}
For  any  integer  $n$  there  is  a  canonical
isomorphism of central extensions
$$[n]^*G\simeq G^{\frac{n(n+1)}{2}}\ot [-1]^*G^{\frac{n(n-1)}{2}}$$
where $[n]^*G$ is the pull-back of $G$  with  respect  to  the
multiplication by $n$ morphism $[n]:K\ra K$. In particular,
if $G$ is symmetric then $[n]^*G\simeq G^{n^2}$.
\end{lem}
      
\Pf  .  The structure of the central extension $G$ of $K$ by
$\G_m$  is  equivalent  to  the  following  data  (see  e.g.
\cite{Breen}): a cube structure on $\G_m$-torsor $G$ over
$K$ and a trivialization of the corresponding biextension
$\La(G)=(p_1+p_2)^*G\ot p_1^*G^{-1}\ot p_2^*G^{-1}$ of $K^2$.
Now for any cube structure there is a canonical isomorphism
(see \cite{Breen})
$$[n]^*G\simeq G^{\frac{n(n+1)}{2}}\ot [-1]^*G^{\frac{n(n-1)}{2}}$$
which is compatible with the natural isomorphism of biextensions
$$([n]\times [n])^*\La(G)\simeq\La(G)^{n^2}\simeq
\La(G)^{\frac{n(n+1)}{2}}\ot
([-1]\times [-1])^*\La(G)^{\frac{n(n-1)}{2}}.$$
The latter isomorphism is compatible with the trivializations
of both sides when $G$ arises from a central extension.
\ed
      
\begin{rem} Locally one can choose a splitting $K\ra  G$  so
that the central extension is given by a 2-cocycle
$f:K\times  K\ra\G_m$.  The previous lemma says that for any
2-cocycle $f$ the functions $f(nk,nk')$ and
$f(k,k')^{\frac{n(n+1)}{2}}f(-k,-k')^{\frac{n(n-1)}{2}}$
differ by a canonical coboundary. In fact this coboundary can
be written explicitly in terms of the functions $f(mk,k)$ for various
$m\in\Z$.
\end{rem}
      
\begin{prop}\label{order}
Assume that $K$ is annihilated by an integer $N$.
If $N$ is odd then for any Heisenberg group $G\ra K$
the central extension $G^N$ is canonically trivial,
otherwise $G^{2N}$ is trivial.
If $G$ is symmetric and $N$ is odd
then $G^N$ (resp. $G^{2N}$ if $N$ is even)
is trivial as a symmetric extension.
\end{prop}
      
\Pf . Combining the previous lemma with (\ref{inver})
we get the following isomorphism:
$$[n]^*G\simeq
G^{\frac{n(n+1)}{2}}\ot (G^{op})^{-\frac{n(n-1)}{2}}
\simeq G^n\ot (G\ot G^{op -1})^{\frac{n(n-1)}{2}}.$$
Now we remark that $G\ot G^{op -1}$ is given by a
trivial $\G_m$-torsor over $K$ with the group law induced
by the commutator form $e:K\times K\ra\G_m$ considered
as 2-cocycle. It remains to note that
$e^{\frac{n(n-1)}{2}}=1$ for $n=2N$ (resp. for $n=N$ if
$N$ is odd). Hence, the triviality of $G^n$ in these cases.
\ed
                      
\begin{cor} Let $G\ra K$ be a symmetric Heisenberg group such
that   the   order   of  $K$  over  $S$  is  odd.  Then  the
$\G_m$-torsor over $K$ underlying $G$ is trivial.
\end{cor}
      
\Pf   .  The  isomorphism  (\ref{inver})  implies  that  the
$\G_m$-torsor over $K$ underlying $G^2$ is trivial. Together
with the  previous  proposition this gives the result.
\ed
       
If  $G\ra  K$ is a (symmetric) Heisenberg group scheme, such
that $K$ is annihilated by an integer $N$,
$n$ is an integer prime to $N$ then $G^n$
is also a (symmetric) Heisenberg group.
When $N$ is odd this group depends only on the residue of $n$ modulo $N$
(due to the triviality of $G^N$).
      
We call a flat subgroup scheme $I\sub   K$ $G$-{\it isotropic} if
the central extension (\ref{ext}) splits over $I$
(in particular, $e|_{I\times I}=1$).
If $\si:I\ra G$ is the corresponding lifting, then
we have the reduced  Heisenberg  group scheme
$$0\ra\G_m\ra p^{-1}(I^\perp)/\si(I)\ra I^\perp/I\ra 0$$
where  $I^\perp\sub  K$  is the orthogonal complement to $I$
with respect to $e$.
If  $G$  is  a symmetric Heisneberg group, then $I\sub K$ is
called {\it symmetrically} $G$-isotropic if the restriction
of  the  central  extension  (\ref{ext})  to  $I$   can   be
trivialized as a symmetric extension. If $\si:I\ra G$  is the
corresponding symmetric lifting them the reduced
Heisenberg group $p^{-1}(I^\perp)/\si(I)$ is also symmetric.
      
Let us define the Witt group $\WH_{\sym}(S)$
as the  group  of  isomorphism  classes  of finite symmetric
Heisenberg  groups  over $S$ modulo the equivalence relation
generated  by  $[G]\sim  [p^{-1}(I^\perp)/\si(I)]$  for  a
symmetrically $G$-isotropic subgroup scheme $I\sub  K$.
The (commutative) addition  in
$\WH_{\sym}(S)$ is defined as follows: if $G_i\ra K_i$ ($i=1,2$) are
Heisenberg groups with commutator forms $e_i$
then their sum is the central extension
$$0\ra \G_m\ra G_1\times_{\G_m} G_2\ra K_1\times K_2\ra 0$$
so that the corresponding commutator form on $K_1\times K_2$
is $e_1\oplus e_2$. The neutral element is the class of $\G_m$
considered as an extension of the trivial group.
The inverse element to $[G]$ is $[G^{-1}]$.
Indeed, there is a canonical splitting of $G\times_{\G_m} G^{-1}\ra
K\times K$ over the diagonal  $K\sub  K\times  K$, hence
the triviality of $[G]+[G^{-1}]$.
We define the order of a finite Heisenberg group scheme $G\ra K$ over $S$
to be the order of $K$ over $S$ (specializing to a  geometric
point  of  $S$  one can see easily that this number has form
$d^2$).
Let us denote by $\WH'_{\sym}(S)$ the analogous Witt group of
finite  Heisenberg group schemes $G$ over $S$ of odd order.
Let also  $\WH(S)$  and  $\WH'(S)$  be  the  analogous  groups
defined   for   all   (not   necessarily  symmetric)  finite
Heisenberg groups over $S$ (with equivalence relation  given
by $G$-isotropic subgroups).

\begin{rem} Let us denote by $\W(S)$ the Witt group of finite
flat group schemes over $S$ with non-degenerate symplectic
$\G_m$-valued forms (modulo the equivalence relation given
by global isotropic flat subgroup schemes).
Let also $\W'(S)$ be the analogous group for group schemes of
odd order. Then we have a natural homomorphism
$\WH(S)\ra \W(S)$  and one
can show that the induced map $\WH'_{\sym}\ra \W'(S)$ is an
isomorphism. This follows essentially from the fact
that a finite symmetric Heisenberg group of odd order is determined
up to an isomorphism by the  corresponding  commutator  form,
also if $G\ra K$ is a  symmetric  finite  Heisenberg
group with the commutator form $e$,
$I\sub K$ is an isotropic flat subgroup scheme of odd order,
then there is a unique symmetric lifting $I\ra G$.
\end{rem}
                                                         
\begin{thm}\label{annih}
The   group   $\WH_{\sym}(S)$  (resp.  $\WH'_{\sym}(S)$)  is
annihilated by $8$ (resp. $4$).
\end{thm}
      
\Pf . Let $G\ra K$ be a symmetric finite Heisenberg group.
Assume first that the order $N$ of $G$ is odd. Then
we can find integers $m$ and $n$ such that
$m^2+n^2\equiv -1\mod(N)$. Let $\a$ be an automorphism of
$K\times K$ given by a matrix
$\left( \matrix m & -n\\ n & m \endmatrix \right)$. Let
$G_1=G\times_{\G_m} G$ be a Heisenberg extension of $K\times
K$ representing the class $2[G]\in \WH'_{\sym}(S)$. Then
from Lemma \ref{mult} and Proposition \ref{order} we get
$\a^*G_1\simeq G_1^{-1}$, hence $2[G]=-2[G]$, i.e. $4[G]=0$
in $\WH'(S)$.
      
If $N$ is even we can apply the similar argument to the
4-th cartesian power of $G$ and the automorphism of $K^4$
given by an integer $4\times 4$-matrix $Z$ such that
$Z^t Z=(2N-1)\id$. Such a matrix can be found by considering the
left multiplication by a quaternion $a+bi+cj+dk$ where
$a^2+b^2+c^2+d^2=2N-1$.
\ed

\section{Schr\"odinger representations}\label{Schr}
      
Let  $G$  be a finite Heisenberg group scheme of order $d^2$
over  $S$. A representation of $G$ of weight 1 is a locally
free $\O_S$-module together with the action of $G$ such that
$\G_m\sub G$ acts by the identity character.
We refer to chapter V of \cite{MB} for basic facts about such
representations. In this section we study the problem of
existence of a {\it Schr\"odinger representation}  for  $G$,
i.~e. a weight-1 representation of $G$ of
rank $d$ (the minimal possible rank).
It is well known that such a representation exists if $S$ is the
spectrum of an  algebraically  closed field (see e.g.
\cite{MB}, V, 2.5.5). Another example is the following.
As we already mentioned one can associate a
finite Heisenberg group scheme $G(L)$
(called the Mumford group) to a line bundle $L$ on an
abelian scheme $\pi:A\ra S$ such that $K(L)$ is finite.
Assume that the base scheme $S$ is connected.
Then $R^i\pi_*(L)=0$ for $i\neq i(L)$ for some integer $i(L)$
(called  the {\it  index}  of  $L$)  and $R^{i(L)}\pi_*(L)$ is a
Schr\"odinger representation for $G(L)$ (this follows from
\cite{Mum} III, 16 and \cite{Muk2}, prop.1.7).
In general, L.  Moret-Bailly  showed in \cite{MB} that a
Schr\"odinger representation
exists after some smooth base change. The main result of this
section is that for symmetric Heisenberg group schemes of odd
order a Schr\"odinger representation always exists.
      
Let  $G$  be  a  symmetric finite Heisenberg group scheme of
order $d^2$ over
$S$.  Then  locally  (in  {\it fppf}  topology) we can choose a
Schr\"odinger representation $V$ of $G$.
According to Theorem V, 2.4.2 of \cite{MB}
for any  weight-1 representation $W$ of $G$ there is
a canonical isomorphism $V\ot\underline{\Hom}_G(V,W)\wt{\ra}
W$. In particular, locally $V$ is unique up to an isomorphism and
$\underline{\Hom}_G(V,V)\simeq\O$. Choose an open covering
$U_i$ such that there  exist  Schr\"odinger  representations
$V_i$ for $G$ over $U_i$. For a sufficently fine  covering
we   have   $G$-isomorphisms   $\phi_{ij}:V_i\ra   V_j$   on
the intersections $U_i\cap U_j$, and
$\phi_{jk}\phi_{ij}=\a_{ijk}\phi_{ik}$   on    the    triple
intersections $U_i\cap U_j\cap U_k$ for some functions
$\a_{ijk}\in\O^*(U_i\cap U_j\cap U_k)$. Then $(\a_{ijk})$ is a
Cech 2-cocycle with values in $\G_m$ whose cohomology class
$e(G)\in H^2(S,\G_m)$ doesn't depend on the choices made.
Furthermore,  by definition $e(G)$ is trivial if and only if
there exists a global weight-1 representation we are  looking
for.
      
Using the language of gerbs (see e.g. \cite{Gir})
we can rephrase the construction above  without  fixing an open
covering. Namely, to each finite Heisenberg group $G$ we can
associate  the  $\G_m$-gerb  $\Sch_G$  on  $S$   such   that
$\Sch_G(U)$  for  an  open  set $U\sub S$ is the category of
Schr\"odinger representations for $G$ over $U$. Then $\Sch_G$
represents the cohomology class $e(G)\in H^2(S,\G_m)$.
      
Notice that the class $e(G)$ is actually represented by
an Azumaya algebra $\AA(G)$ which is defined as follows.
Locally,  we  can  choose a Schr\'odinger representation $V$
for $G$ and put $\AA(G)=\underline{\End}(V)$.  Now  for  two
such  representations  $V$  and  $V'$  there  is a canonical
isomorphism of algebras
$\underline{\End}(V)\simeq\underline{\End}(V')$ induced
by  any  $G$-isomorphism  $f:V\ra   V'$   (since   any  other
$G$-isomorphism differs from $f$ by a scalar), hence these
local  algebras glue together into a global Azumaya algebra
$\AA(G)$ of rank $d^2$. In particular, $d\cdot e(G)=0$  (see  e.g.
\cite{Groth1},  prop. 1.4).
      
Now let $W$ be a {\it global}
weight-1 representation of $G$ which is locally free of rank
$l\cdot d$ over $S$. Then we claim that $\underline{\End}_G(W)$
is an Azumaya algebra with the class  $-e(G)$.  Indeed,
locally  we  can  choose a representation $V$ of rank $d$ as
above and a $G$-isomorphism $W\simeq V^l$ which induces a
local isomorphism $\underline{\End}_G(W)\simeq\Mat_l(\O)$.
Now we claim that there is a global algebra isomorphism
$$\AA(G)\ot\underline{\End}_G(W)\simeq\underline{\End}(W).$$
Indeed, we have
canonical isomorphism of $G$-modules of  weight  1
(resp. $-1$)
$V\ot\und{\Hom}_G(V,W)\wt{\ra}W$  (resp.
$V^*\ot\und{\Hom}_G(V^*,W^*)\wt{\ra}W^*$). Hence,
we have a sequence of natural morphisms
\begin{align*}
&\und{\End}(W)\simeq W^*\ot W\simeq
V^*\ot V\ot\und{\Hom}_G(V^*,W^*)\ot\und{\Hom}_G(V,W)\ra\\
&\ra\und{\End}(V)\ot\Hom_{G\times G}(V^*\ot V,W^*\ot W)\ra
\und{\End}(V)\ot\und{\End}_G(W)
\end{align*}
--- the latter map is obtained by taking the image of the
identity   section   $\id\in  V^*\ot  V$  under  a  $G\times
G$-morphism $V^*\ot V\ra W^*\ot W$. It is easy to see that
the composition morphism gives the required isomorphism.
This leads to the following statement.
      
\begin{prop}\label{obst}
For any finite Heisenberg group scheme $G$ over $S$
a canonical element $e(G)\in\Br(S)$  is  defined  such  that
$e(G)$   is   trivial   if   and  only  if  a  Schr\"odinger
representation for $G$ exists. Furthermore, $d\cdot e(G)=0$
where the order of $G$ is $d^2$, and if there  exists
a  weight-1 $G$ representation which is locally free of rank
$l\cdot d$ over $S$ then $l\cdot e(G)=0$.
\end{prop}
                
\begin{prop}\label{hom}
The map $[G]\mapsto e(G)$ defines a homomorphism
$\WH(S)\ra\Br(S)$.
\end{prop}
      
\Pf  .  First  we  have  to  check  that  if  $I\sub K$ is a
$G$-isotropic subgroup, $\wt{I}\sub G$ its lifting, and
$\ov{G}=p^{-1}(I^{\perp})/\wt{I}$ then $e(\ov{G})=e(G)$.
Indeed, there is a canonical equivalence of $\G_m$-gerbs
$\Sch_G\ra\Sch_{\ov{G}}$ given by the functor
$V\mapsto V^{\wt{I}}$ where $V$ is a  (local)  Schr\"odinger
representation of $G$. Next if $G=G_1\times_{\G_m} G_2$,
then   for   every   pair   $(V_1,V_2)$   of weight-1
representations of $G_1$ and $G_2$ there is a natural
structure of weight-1  $G$-representation on $V_1\ot V_2$,
hence we get an equivalence of $\G_m$-gerbs
$\Sch_{G_1}+\Sch_{G_2}\ra\Sch_G$ which implies the equality
$e(G)=e(G_1)+e(G_2)$. At last, the map
$V\ra        V^*$        induces        an       equivalence
$\Sch_{G}^{op}\ra\Sch_{G^{-1}}$ so that $e(G^{-1})=-e(G)$.
\ed
         
\begin{thm}\label{odd}
Let $G$ be a symmetric finite  Heisenberg  group scheme
of  odd order.  Then $e(G)=0$, that is there exists a global
Schr\"odinger representation for $G$.
\end{thm}
      
\Pf . Let $[G]\in\WH'_{\sym}(S)$ be a class of $G$ in the Witt group.
Then $4[G]=0$ by Theorem \ref{annih}, hence $4e(G)=0$ by
Proposition \ref{hom}. On the other hand, $d\cdot e(G)=0$ by
Proposition \ref{obst} where $d$ is odd, therefore, $e(G)=0$.
\ed
      
Let us give an example of a symmetric finite Heisenberg group
scheme  of  {\it  even}  order   without   a   Schr\"odinger
representation. First let us recall the construction from
\cite{sympl}  which  associates to a  group
scheme $G$ over $S$ which is a central extension of a finite commutative
group  scheme  $K$ by $\G_m$, and a $K$-torsor $E$ over $S$ a
class $e(G,E)\in H^2(S,\G_m)$. Morally, the map
$$H^1(S,K)\ra H^2(E,\G_m): E\mapsto e(G,E)$$
is the boundary homomorphism corresponding
to the exact sequence
$$0\ra \G_m\ra G\ra K\ra 0.$$
To define it consider the category $\CC$ of liftings of $E$ to
to a $G$-torsor. Locally such a lifting always exists and any
two such liftings differ by a $\G_m$-torsor. Thus, $\CC$
is a $\G_2$-gerb over $S$, and by definition
$e(G,E)$ is the class of $\CC$ in $H^2(S,\G_m)$
Note that $e(G,E)=0$ if and only if there exists a $G$-equivariant
line bundle $L$ over $E$, such that $\G_m\sub G$ acts on $L$
via the identity character.

A $K$-torsor $E$
defines a commutative group extension $G_E$ of $K$ by $\G_m$
as follows. Choose local trivializations of  $E$  over  some
covering  $(U_i)$ and let $\a_{ij}\in K(U_i\cap U_j)$ be the
corresponding 1-cocycle with values in $K$. Now we glue $G_E$
from the trivial extensions $\G_m\times K$ over $U_i$ by the
following transition isomorphisms over $U_i\cap U_j$:
$$f_{ij}:\G_m\times K\ra\G_m\times K:(\la,x)\mapsto
(\la e(x,\a_{ij}),x)$$
where  $e:K\times   K\ra\G_m$   is   the   commutator   form
corresponding  to  $G$. It is easy to see that $G_E$ doesn't
depend on a choice of trivializations. Now we claim that
if $G$ is a Heisenberg group then
\begin{equation}\label{diff}
e(G,E)=e(G\ot G_E)-e(G).
\end{equation}
This is checked by a direct computation with Cech cocycles.
Notice that if $E^2$ is a trivial $K$-torsor
then $G_E^2$ is a trivial central extension  of  $K$,  hence
$G_E$ is a symmetric extension. Thus, if $G$ is a symmetric
Heisenberg group, then $G\ot G_E$ is also symmetric.
As was shown in \cite{sympl} the left hand side of (\ref{diff})
can be non-trivial.
Namely, consider the case when $S=A$ is a principally
polarized abelian variety over an algebraically closed field
$k$ of characteristic $\neq2$. Let $K=A_2\times A$ considered
as a (constant) finite group scheme over $A$. Then we can
consider $E=A$ as a $K$-torsor over $A$ via the morphism
$[2]:A\ra A$. Now if $G\ra A_2$ is a Heisenberg extension of
$A_2$ (defined over $k$) then we can consider $G$ as
a constant group scheme over $A$ and the class $e(G,E)$
is trivial if and only if $G$ embeds into the Mumford
group $G(L)$ of some line bundle $L$ over $A$ (this
embedding should be the identity on $\G_m$).
When $\NS(A)=\Z$ this means, in particular, that
the commutator form $A_2\times A_2\ra\G_m$
induced by $G$ is proportional to the symplectic
form given by the principal polarization.
When $\dim A\ge 2$ there is a plenty of other
symplectic forms on $A_2$, hence, $e(G,E)$ can be non-trivial.
                                                              
Now we are going to show that one can replace $A$ by its
general point in this example. In other words, we
consider the base $S=\Spec(k(A))$ where $k(A)$ is the field
of rational functions on $A$. Then $E$ gets replaced by
$\Spec(k(A))$ considered as a $A_2$-torsor over itself
corresponding to the Galois extension
$$[2]^*:k(A)\ra k(A): f\mapsto f(2\cdot)$$
with the Galois group $A_2$. Note that the class $e(G,E)$
for any Heisenberg extension $G$ of $A_2$ by $k^*$
is annihilated by the pull-back to $E$, hence, it is
represented by the class of Galois cohomology
$H^2(A_2,k(A)^*)\sub\Br(k(A))$ where $A_2$ acts
on $k(A)$ by translation of argument. It is easy
to see that this class is the image of the class
$e_G\in H^2(A_2,k^*)$ of the central extension $G$
under the natural homomorphism
$H^2(A_2,k^*)\ra H^2(A_2,k(A)^*)$.
From the exact sequence of groups
$$0\ra k^*\ra k(A)^*\ra k(A)^*/k^*\ra 0$$
we get the exact sequence of cohomologies
$$0\ra H^1(A_2,k(A)^*/k^*)\ra H^2(A_2,k^*)\ra
H^2(A_2,k(A)^*)$$
(note that $H^1(A_2,k(A)^*)=0$ by Hilbert theorem 90).
It follows that central extensions $G$ of $A_2$ by
$k^*$ with trivial $e(G,E)$ are classified by elements
of $H^1(A_2,k(A)^*/k^*)$.

\begin{lem} Let $A$ be a principally polarized abelian variety
over an algebraically closed field $k$  of characteristic
$\neq 2$. Assume that $\NS(A)=\Z$. Then
$H^1(A_2,k(A)^*/k^*)=\Z/2\Z$.
\end{lem}

\Pf . Interpreting $k(A)^*/k^*$ as the group of divisors
linearly equivalent to zero we obtain the exact sequence
$$0\ra k(A)^*/k^*\ra\Div(A)\ra\Pic(A)\ra 0,$$
where $\Div(A)$ is the group of all divisors on $A$.
Note that as $A_2$-module $\Div(A)$ is decomposed into
a direct sum of modules of the form $\Z^{A_2/H}$ where
$H\sub A_2$ is a subgroup. Now by Shapiro lemma we have
$H^1(A_2,\Z^{A_2/H})\simeq H^1(H,\Z)$, and the latter group
is zero since $H$ is a torsion group. Hence,
$H^1(A_2,\Div(A))=0$. Thus, from the above exact sequence
we get the identification
$$H^1(A_2,k(A)^*/k^*)\simeq
\coker(\Div(A)^{A_2}\ra\Pic(A)^{A_2}).$$
Now we use the exact sequence
$$0\ra\Pic^0(A)\ra\Pic(A)\ra\NS(A)\ra 0,$$
where $\Pic^0(A)=\hat{A}(k)$.
Since the actions of $A_2$ on $\Pic^0(A)$ and $\NS(A)$ are trivial
we have the induced exact sequence
$$0\ra\Pic^0(A)\ra\Pic(A)^{A_2}\ra\NS(A).$$
The image of the right arrow is the subgroup $2\NS(A)\sub\NS(A)$.
Note that $\Pic^0(A)=[2]^*\Pic^0(A)$, hence this subgroup
belongs to the image of $[2]^*\Div(A)\sub\Div(A)^{A_2}$.
Thus, we deduce that
$$H^1(A_2,k(A)^*/k^*)\simeq\coker(\Div(A)^{A_2}\ra 2\NS(A)).$$
Let  $[L]\sub\NS(A)$  be  the  generator corresponding to
a line bundle $L$ of degree 1 on $A$.
Then $L^4=[2]^*L$, hence $4\cdot [L]=[L^4]$ belongs to the image
of $\Div(A)^{A_2}$. On the other hand, it is easy to see
that there is no $A_2$-invariant divisor representing
$[L^2]$, hence
$$H^1(A_2,k(A)^*/k^*)\simeq\Z/2\Z.$$
\ed

It follows that under the assumptions of this lemma there
is a unique Heisenberg extensions $G$ of $A_2$ by $k^*$
with the trivial class $e(G,E)$ (the Mumford
extension corresponding to $L^2$,
where $L$ is a line bundle of degree 1 on $A$). Hence,
for $g\ge 2$ there exists a Heisenberg extension with
a non-trivial class $e(G,E)\in\Br(k(A))$.

\section{Representations of the Heisenberg groupoid}
                                             
Recall that the Heisenberg group $H(W)$ associated with
a symplectic vector space $W$ is a central extension
$$0\ra T\ra H(W)\ra W\ra 0$$
of $W$ by the 1-dimensional torus $T$ with the commutator form
$\exp(B(\cdot,\cdot))$ where $B$ is the symplectic form.
In this section we consider an analogue of this
extension  in  the   context   of   abelian   schemes   (see
\cite{Weilrep} , sect. 7, \cite{sympl}).  Namely,  we
replace a vector space $W$ by an abelian scheme $X/S$. Bilinear
forms on $W$ get replaced by biextensions of $X^2$.
Recall that
a {\it biextension} of $X^2$ is a line bundle $\LL$ on $X^2$
together with isomorphisms
\begin{align*}
&\LL_{x+x',y}\simeq \LL_{x,y}\ot \LL_{x',y},\\
&\LL_{x,y+y'}\simeq \LL_{x,y}\ot \LL_{x,y'}
\end{align*}
--- this is a symbolic notation for isomorphisms
$(p_1+p_2,p_3)^*\LL\simeq p_{13}^*\LL\ot p_{23}^*\LL$ and
$(p_1,p_2+p_3)^*\LL\simeq p_{12}^*\LL\ot p_{13}^*\LL$ on $X^3$,
satisfying some natural cocycle conditions (see e.g. \cite{Breen}).
The parallel notion to the skew-symmetric form on $W$ is that of
a {\it  skew-symmetric  biextension}  of  $X^2$ which is  a
biextension $\LL$ of $X^2$ together with
an   isomorphism   of   biextensions
$\phi:\si^*\LL\wt{\ra}  \LL^{-1}$, where   $\si:X^2\ra  X^2$  is  the
permutation of factors, and a trivialization
$\De^*\LL\simeq\O_X$ of $\LL$ over the diagonal
$\De:X\ra X^2$ compatible with $\phi$.
A skew-symmetric biextension $\LL$ is called {\it symplectic}
if the corresponding homomorphism $\psi_{\LL}:X\ra\hat{X}$ (where
$\hat{X}$ is the dual abelian scheme) is an isomorphism.
An {\it isotropic} subscheme (with respect to $\LL$) is an abelian
subscheme $Y\sub X$ such that there is an isomorphism of skew-symmetric
biextensions  $\LL|_{Y\times Y}\simeq\O_{Y\times Y}$.
This is equivalent to the condition that the composition
$Y\stackrel{i}{\ra} X\stackrel{\psi_{\LL}}{\ra}\hat{X}
\stackrel{\hat{i}}{\ra} \hat{Y}$ is zero.
An isotropic subscheme $Y\sub X$ is called {\it lagrangian}
if the morphism $Y\ra \ker(\hat{i})$ induced by $\psi_{\LL}$ is an
isomorphism. In particular, for such a subscheme the quotient
$X/Y$ exists and is isomorphic to $\hat{Y}$.
      
Note that to define the Heisenberg group extension it is not sufficient
to have a symplectic form $B$ on $W$: one needs a bilinear form $B_1$
such that $B(x,y)=B_1(x,y)-B_1(y,x)$. In the case of the real symplectic
space one can just take $B_1=B/2$, however in our situation we have
to simply add necessary data.
An {\it enhanced} symplectic biextension $(X,\BB)$  is a
biextension $\BB$ of $X^2$ such that
$\LL:=\BB\ot\si^*\BB^{-1}$ is a symplectic biextension.
The standard enhanced symplectic biextension
for $X=\hat{A}\times A$, where $A$ is any abelian scheme, is
obtained by setting
$$\BB=p_{14}^*\PP\in\Pic(\hat{A}\times A\times \hat{A}\times A),$$
where $\PP$ is the normalized Poincar\'e line bundle on $A\times\hat{A}$.
     
Given an enhanced symplectic
biextension   $(X,\BB)$  one  defines  the  {\it  Heisenberg
groupoid} $H(X)=H(X,\BB)$ as
the stack of monoidal groupoids
such that $H(X)(S')$ for an $S$-scheme $S'$
is the monoidal groupoid generated by
the central subgroupoid $\PPic(S')$ of $\G_m$-torsors on $S'$ and the symbols
$T_x$, $x\in X(S')$ with the composition law
$$T_x\circ T_{x'}= \BB_{x,x'} T_{x+x'}.$$
The Heisenberg groupoid is a central extension of $X$ by the
stack of line bundles on $S$ in the sense of Deligne \cite{Des}.
    
In \cite{Weilrep} we considered the action of $H(\hat{A}\times A)$
on $\D^b(A)$ which is similar to the standard representation of
the Heisenberg group $H(W)$ on functions on a lagrangian subspace of
$W$. Below we construct similar representations of the Heisenberg
groupoid $H(X)$ associated with lagrangian subschemes in $X$.
Further,  we  construct  intertwining  functors for two such
representations
corresponding to a pair of lagrangian subschemes, and consider
the analogue of Maslov index for a triple of lagrangian subschemes
that arises when composing these intertwining functors.
    
To define an action of $H(X)$ associated with a lagrangian subscheme one
needs some  auxilary data described as follows.
An  {\it  enhanced}  lagrangian  subscheme  (with respect to
$\BB$) is a pair $(Y,\a)$ where $Y\sub X$ is  a  lagrangian
subscheme with respect to $X$, $\a$ is a line bundle on $Y$
with a rigidification along the zero section such that an
isomorphism of symmetric
biextensions $\La(\a)\simeq \BB|_{Y\times Y}$ is given, where
$\La(\a)=(p_1+p_2)^*\a\ot p_1^*\a^{-1}\ot p_2^*\a^{-1}$.
Note that an enhanced lagrangian subscheme is a particular
case of an {\it isotropic pair} as defined in
\cite{Weilrep} II, 7.3.
                            
With every enhanced lagrangian subscheme $(Y,\a)$
one can associate a representation of $H(X)(S)$  as  follows
(see \cite{Weilrep},\cite{sympl}). Let
$\DD(Y,\a)$ be the category of pairs $(\FF,a)$ where
$\FF\in\D^b(X)$, $a$ is an isomorphism in $\D^b(Y\times X)$:
\begin{equation}\label{Schrsp}
a:(i_Yp_1+p_2)^*\FF\wt{\ra} \BB^{-1}|_{Y\times X}
\ot p_1^*\a^{-1}\ot p_2^*\FF
\end{equation}
where $i_Y:Y\hra X$ is the embedding, such that $(e\times\id)^*a=\id$.
These data should satisfy
the following cocycle condition:
$$(p_1+p_2,p_3)^*a=(p_2,p_3)^*a\circ (p_1,i_Yp_2+p_3)^*a$$
in $\D^b(Y\times Y\times X)$.
Then there is a natural action of  the  Heisenberg  groupoid
$H(X)(S)$
on  the  category  $\DD(Y,\a)$ such that a line bundle $M$ on $S$
acts by tensoring with $p^*M$ and a generator $T_x$ acts
by the functor
\begin{equation}\label{act}
\FF\mapsto \BB|_{X\times x}\ot t_x^*(\FF).
\end{equation}
If  $S'$  is  an  $S$-scheme then this action is compatible
with the action of $H(X)(S')$ on $\DD(Y_{S'},\a_{S'})$ via
pull-back functors.
      
Let  $\de_{Y,\a}\in\DD(Y,\a)$  be   the   following   object
(delta-function at $(Y,\a)$):
\begin{equation}\label{delta}
\de_{Y,\a}=i_{Y*}(\a^{-1})
\end{equation}
where $i_Y:Y\ra X$ is the embedding. It is easy to see that
$\de_{Y,\a}$ has a canonical structure of an object of $\DD(Y,\a)$
and         for        $y\in        Y$        one        has
$T_y(\de_{Y,\a})\simeq\a_y^{-1}\de_{Y,\a}$.
      
Let  $(Y,\a)$,  $(Z,\b)$  be  a  pair of enhanced lagrangian
subschemes in $X$, such that $Y\cap Z$ is finite over $S$. Then
the natural morphism $Y\ra X/Z\simeq\hat{Z}$ is an  isogeny,
hence, $Y\cap Z$ is flat over $S$. Note that we have isomorphisms of
biextensions $\La(\a|_{Y\cap Z})\simeq\La(\b|_{Y\cap Z})\simeq
\BB|_{(Y\cap Z)^2}$, hence the trivialization of
$\La(\b|_{Y\cap Z}\ot \a^{-1}|_{Y\cap Z})$. Thus, the
$\G_m$-torsor $G_{Y,Z}=\b|_{Y\cap Z}\ot \a^{-1}|_{Y\cap Z}$
has a natural structure of a central extension of
$Y\cap Z$ by $\G_m$. Furthermore, the corresponding commutator
form   $(Y\cap   Z)^2\ra\G_m$  is  non-degenerate  since  it
corresponds to the canonical duality between
$Y\cap Z=\ker(Y\ra\hat{Z})$ and $Y\cap Z=\ker(Z\ra\hat{Y})$
(see \cite{sympl}, remark after Prop. 3.1).
Thus, $G_{Y,Z}$ is a finite Heisenberg group scheme over $S$.
If the line bundles $\a$ and $\b$
are symmetric then so is $G_{Y,Z}$.
      
Let $V$ be a Schr\"odinger representation of $G_{Y,Z}$.
Generalizing the construction of \cite{sympl}
we define the $H(X)(S)$-intertwining operator
$$R(V):\DD(Y,\a)\ra\DD(Z,\b):
\FF\mapsto
\und{\Hom}_{G_{Y,Z}}(V,p_{2*}(\BB|_{Z\times X}\ot p_1^*\b\ot
(i_Zp_1+p_2)^*\FF)).$$
Here $p_1$ and $p_2$ are the projections of the product
$Z\times_S X$ onto its factors. The $G_{Y,Z}$-module structure
on   $p_{2*}(\BB|_{Z\times  X}\ot  p_1^*\b\ot(i_Zp_1+p_2)^*\FF)$
comes from the natural $G_{Y,Z}$-action on
$I(\FF)=\BB|_{Z\times  X}\ot  p_1^*\b\ot(i_Zp_1+p_2)^*\FF$ which is
compatible with the action of $Y\cap Z$ on  $Z\times  X$  by
the translation  of  the  first  argument  and  arises  from  the
canonical isomorphism
\begin{equation}\label{integrand}
I(\FF)_{(z+u,x)}\simeq \b_u\a_u^{-1} I(\FF)_{(z,x)}
\end{equation}
where $z\in  Z$,  $x\in  X$,  $u\in  Y\cap  Z$  (one  should
consider  this  as  an  isomorphism in $\D^b((Y\cap Z)\times
Z\times X)$). When $V$ is the representation
associated with a lagrangian subgroup scheme $H\sub G_{Y,Z}$
this functor coincides with the one defined in \cite{sympl}.
      
Let us call an enhanced lagrangian subscheme $(Y,\a)$
{\it admissible} if the projection $X\ra X/Y$ splits.
For such a subscheme we have an equivalence
$\DD(Y,\a)\simeq\D^b(X/Y)$. Namely, let $s_{X/Y}:X/Y\ra X$ be a
splitting of the canonical projection $q_{X/Y}:X\ra X/Y$. Let
$q_Y=\id-s_{X/Y}q_{X/Y}:X\ra Y$ be the corresponding projection to $Y$.
Then the functors $\FF\mapsto s_{X/Y}^*\FF$ and $\GG\mapsto
(q_Y,s_{X/Y}q_{X/Y})^*\BB^{-1}\ot q_Y^*\a^{-1}\ot q_{X/Y}^*\GG$    where
$\FF\in\DD(Y,\a)$,  $\GG\in\D^b(X/Y)$  give   the   required
equivalence.
When  $(Y,\a)$  and  $(Z,\b)$  are  both  admissible  we can
represent the above functor $R(V):\D^b(X/Y)\ra\D^b(X/Z)$
in the standard "integral" form.
      
\begin{lem}\label{compker}
Assume   that   $(Y,\a)$   and   $(Z,\b)$   are
admissible, $Y\cap Z$ is finite. Then
$R(V)(\GG)\simeq p_{2*}(p_1^*\GG\ot \KK(V))$ where
$p_i$ are the projections of $X/Y\times X/Z$ on its factors,
$\KK(V)$ is the following vector bundle on $X/Y\times X/Z$:
\begin{eqnarray}\label{kernel}
\KK(V)=(p_1-q_{X/Y}s_{X/Z}p_2)^*E(V)\ot
(s_{X/Y}p_1-s_{X/Z}p_2,s_{X/Y}p_2)^*\BB\ot \nonumber\\
(s_{X/Y}(p_1-q_{X/Y}s_{X/Z}p_2),q_Ys_{X/Z}p_2)^*\BB\ot
(q_Ys_{X/Z}p_2)^*\a^{-1}
\end{eqnarray}
where  $s_{X/Y}:X/Y\ra  X$ (resp. $s_{X/Z}:X/Z\ra X$) is the
splitting  of  the  projection  $q_{X/Y}:X\ra  X/Y$   (resp.
$q_{X/Z}:X\ra X/Z$), $q_Y=\id-s_{X/Y}q_{X/Y}$, $E(V)$ is the
following bundle on $X/Y$:
\begin{equation}\label{kernelaux}
E(V)=\und{\Hom}_{G_{Y,Z}}(V,(q_{X/Y}i_Z)_*(\b\ot (q_Yi_Z)^*\a^{-1}\ot
(i_Z,s_{X/Y}q_{X/Y}i_Z)^*\BB^{-1}))
\end{equation}
where $i_Y:Y\ra X$, $i_Z:Z\ra X$ are the embeddings.
\end{lem}
      
\Pf . By definition we have
\begin{align*}
&R(V)(\GG)_{\bar{x}}\simeq\\
&\und{\Hom}(V,\int_Z \b_z \BB_{z,s_{X/Z}(\bar{x})}
\BB^{-1}_{q_Y(z+s_{X/Z}(\bar{x})),s_{X/Y}q_{X/Y}(z+s_{X/Z}
(\bar{x}))}\a^{-1}_{q_Y(z+s_{X/Z}(\bar{x}))}\GG_{q_{X/Y}(z+s_{X/Z}
(\bar{x}))} dz)
\end{align*}
where $\bar{x}\in X/Z$, $z\in Z$. Using the isomorphism
$\a_{q_Y(z+s_{X/Z}(\bar{x}))}\simeq
\a_{q_Y(z)}\a_{q_Ys_{X/Z}(\bar{x})}
\BB_{q_Y(z),q_Ys_{X/Z}(\bar{x})}$ and collecting together terms
depending only on $\bar{z}=q_{X/Y}(z)$ we get
\begin{align*}
&R(V)(\GG)_{\bar{x}}\simeq
\und{\Hom}(V,\int_Z   \b_z    \a_{q_Y(z)}^{-1}
\BB^{-1}_{z, s_{X/Y}\bar{z}}    \\
&\Bigl(
\BB_{s_{X/Y}(\bar{z}+q_{X/Y}s_{X/Z}(\bar{x}))-s_{X/Z}(\bar{x}),
s_{X/Y}(\bar{z}+q_{X/Y}s_{X/Z}(\bar{x}))}
\BB_{s_{X/Y}(\bar{z}),q_Ys_{X/Z}(\bar{x})}
\a^{-1}_{q_Ys_{X/Z}(\bar{x})}
\GG_{\bar{z}+q_{X/Y}s_{X/Z}(\bar{x})}\Bigr) dz)\simeq\\
&\int_{X/Y} E(V)_{\bar{z}}
\BB_{s_{X/Y}(\bar{z}+q_{X/Y}s_{X/Z}(\bar{x}))-s_{X/Z}(\bar{x}),
s_{X/Y}(\bar{z}+q_{X/Y}s_{X/Z}(\bar{x}))}
\BB_{s_{X/Y}(\bar{z}),q_Ys_{X/Z}(\bar{x})}\\
&\a^{-1}_{q_Ys_{X/Z}(\bar{x})}
\GG_{\bar{z}+q_{X/Y}s_{X/Z}(\bar{x})} d\bar{z}
\end{align*}
where $\bar{z}$ is now considered as a variable on $X/Y$.
Making the change of variables $\bar{z}\mapsto
\bar{z}-q_{X/Y}s_{X/Z}(\bar{x})$ we arrive to the formula
(\ref{kernel}).
\ed
      
\begin{thm}\label{eq}
Assume that $(Y,\a)$ and $(Z,\b)$ are admissible
and $Y\cap Z$ is finite, then $R(V)$ is an equivalence of categories.
Let  $(T,\ga)$  be an admissible  enhanced  lagrangian
subscheme such that $Y\cap T$ and $Z\cap T$ are finite,
$W$  (resp.  $U$)  be  a  Schr\"odinger  representation  for
$G_{Y,T}$ (resp. $G_{Z,T}$). Then
$$R(U)\circ R(V)\simeq R(W)\ot M[n]$$
for some line bundle $M$ on $S$ and some integer $n$.
\end{thm}
      
\Pf . The direct computation shows that the
kernel $\KK(V)\in\D^b(X/Y\times X/Z)$ constructed above
satisfies  the "uniform"  intertwining  property
(with respect to $H(X)$-action) defined
in \cite{Weilrep}. Hence,  the
analogue  of  Schur  lemma  for  the  action  of  $H(X)$  on
$\D^b(X/Y)$ where $Y$ is an admissible lagrangian subscheme
(see \cite{Weilrep} Thm 7.9) implies that
$$p_{13*}(p_{12}^*\KK(V)\ot p_{23}^*\KK(V^*))\simeq\De_*(F)$$
where $\KK(V^*)$ is the similar kernel on $X/Z\times X/Y$
giving rise to the functor $R(V^*):\D^b(X/Z)\ra\D^b(X/Y)$,
$p_{ij}$ are the projections of $X/Y\times X/Z\times X/Y$ on
the pairwise products, $\De:X/Y\ra (X/Y)^2$ is the  diagonal
embedding.
In the case when $S$ is the spectrum
of a field we know that $F\simeq N[n]$ for some line bundle
$N$  on  $S$  and  some  integer $n$ (see \cite{Weilrep}). By
Prop.1.7 of \cite{Muk2} this implies that the same  is  true
when $S$ is connected. Therefore, in this case the composition
$R(V^*)\circ R(V)$ is isomorphic to the tensoring with $N[n]$.
Repeating this for the composition $R(V)\circ
R(V^*)$ we conclude that $R(V)$ is an equivalence.
Similar argument works for the proof of the second assertion.
\ed
       
\begin{rems} 1. Most probably, one can extend this theorem
to the case of arbitrary enhanced lagrangian subschemes.
However, it seems that the definition of $\DD(Y,\a)$ should be
modified in this case (one should start with appropriate
category of complexes and then localize  it).
   
\noindent 2. An integer $n$ and a line bundle $M$ on $S$ appearing
in the above theorem should be considered as analogues of
the Maslov index (see \cite{LV}) for a triple  $(Y,Z,T)$.
Note that different choices of Schr\"odinger representations  $V$,
$W$, and $U$
above affect $M$ but not $n$, hence the function $n(Y,Z,T)$
behaves  very  much  like  the  classical  Maslov index
(cf. \cite{Orlov}).
\end{rems}

Let $(Y,\a)$, $(Z,\b)$ and $(T,\ga)$ be a triple of enhanced
lagrangian subschemes in $X$. Let us denote by $K=K(Y,Z,T)$
the  kernel  of  the  homomorphism  $Y\times  Z\times   T\ra
X:(y,z,t)\mapsto y+z+t$. Let $p_Y:K\ra Y$, $p_Z:K\ra Z$ and
$p_T:K\ra T$ be the restrictions to $K$ of the natural projections
from  $Y\times  Z\times  T$  to  its  factors.  Consider the
following line bundle on $K$:
\begin{equation}\label{MYZT}
M(Y,Z,T)=(-p_Y)^*\a^{-1}\ot p_Z^*\b\ot p_T^*\ga\ot
(p_Z,p_T)^*\BB|_{Z\times T}.
\end{equation}
Then  $M(Y,Z,T)$  has  a canonical cube structure induced by
that of $\a$, $\b$, $\ga$ and $\BB$.
                 
\begin{lem} There are canonical isomorphisms of line bundles
with cube structures on $K$
$$M(Y,Z,T)\simeq M(Z,T,Y)\simeq M(T,Y,Z).$$
There is a canonical isomorphism of biextensions of $K\times
K$:
$$\La(M(Y,Z,T))\simeq (p_Zp_2,p_Tp_1)^*\LL$$
where $p_i$ are the projections of $K\times K$ on its factors.
\end{lem}
      
\Pf . We have
$$(M(Z,T,Y)\ot M(Y,Z,T)^{-1})_{y,z,t}=
\a_y\a_{-y}\b_z^{-1}\b_{-z}^{-1}\BB_{t,y}\BB_{z,t}^{-1}\simeq
\BB_{y,y}\BB_{z,z}^{-1}\BB_{t,y}\BB_{z,-t}$$
where    $y+z+t=0$    (here    we   used   the   isomorphism
$\a_y\a_{-y}\simeq  \BB_{y,-y}\simeq  \BB_{y,y}^{-1}$  and   the
similar isomorphism for $\b$).
It is easy to see that when we substitute $t=-y-z$ the right
hand side becomes trivial.
      
The second isomorphism is obtained as follows:
$$\La(M(Y,Z,T))_{(y,z,t),(y',z',t')}\simeq
\BB_{-y,-y'}^{-1}\BB_{z,z'}\BB_{t,t'}\BB_{z,t'}\BB_{z',t}.$$
If we substitute $-y=z+t$, $-y'=z'+t'$ the right hand side
becomes $\BB_{t,z'}^{-1}\BB_{z',t}\simeq \LL_{z',t}$.
\ed
      
Consider the embedding $Z\cap T\hookrightarrow K:u\mapsto
(0,-u,u)$.   Then   the   previous   lemma   implies    that
$\La(M(Y,Z,T))$ is trivial over $(Z\cap T)\times K$. Hence,
$M(Y,Z,T)|_{Z\cap T}$ has a structure of central extension
and the action of $Z\cap T$ on $K$ by translations lifts to
an action of this central extension on $M(Y,Z,T)$.
Moreover,   we  have  a  canonical  isomorphism  of  central
extensions
$$M(Y,Z,T)_{(0,-u,u)}\simeq\ga_u\b_{-u}\BB_{-u,u}\simeq
\ga_u\b_u^{-1}=(G_{Z,T})_u.$$
Hence,  there  is  an  action  of  $G_{Z,T}$  on  $M(Y,Z,T)$
compatible   with   the  action  of  $Z\cap  T$  on  $K$  by
translations. Using cyclic permutation we get embeddings of
$Y\cap T$ and $Y\cap Z$ into $K$ and it is easy to see that
the images of the three embeddings are independent so that
we get an embedding
$(Y\cap Z)\times_S (Z\cap T)\times_S (Y\cap T)\hookrightarrow K$
and the compatible
action of $G_{Y,Z}\times_{\G_m}G_{Z,T}\times_{\G_m} G_{T,Y}$
on $M(Y,Z,T)$.
      
\begin{thm}\label{intert}
With the notation  and  assumptions  of  Theorem
\ref{eq} we have
$$M[n]\simeq\und{\Hom}_{G_{Y,Z,T}}(V_{Y,Z,T},p_*M(Y,Z,T))$$
where $G_{Y,Z,T}=G_{Y,Z}\times_{\G_m} G_{Z,T}\times_{\G_m}
\times G_{T,Y}$, $V_{Y,Z,T}=V\ot U\ot W^*$, $p:K\ra S$ is the
projection.
\end{thm}
      
\Pf   .  Let  us  compare  the  restrictions  of  $R(U)\circ
R(V)(\de)$ and $R(W)(\de)$ to the zero section, where
$\de=\de_{Y,\a}\in\DD(Y,\a)$   is   the  delta-function  at  $Y$
defined by (\ref{delta}). On the one hand, we have
\begin{align*}
&R(U)\circ R(V)(\de)_0\simeq\und{\Hom}_{G_{Y,Z}\times_{\G_m} G_{Z,T}}
(V\ot U,\int_{Z\times T} \ga_t \BB_{z,t}\b_z\de_{z+t} dz dt)\simeq\\
&\simeq\und{\Hom}_{G_{Y,Z}\times_{\G_m} G_{Z,T}}(V\ot U,
\int_{K}M(Y,Z,T)).
\end{align*}
On the other hand,
$$R(W)(\de)_0\simeq\und{\Hom}_{G_{Y,T}}(W,
\int_{Y\cap T} \ga_u\a_u^{-1} du)\simeq W^*$$
since $\int_{Y\cap T}G_{Y,T}\simeq W^*\ot W$
by \cite{MB} V 2.4.2. Therefore,
$$\int_{K} M(Y,Z,T)\simeq V\ot U\ot W^*\ot M[n]$$
as a representation of $G_{Y,Z,T}$.
\ed
      
Consider the following example.
Let $X=\hat{A}\times A$, $\BB=p_{14}^*\PP$,
$(Y,\a)=(A,\O_A)$, $(T,\ga)=(\hat{A},\O_{\hat{A}})$, and
$(Z,\b)=(Z_{\phi,m},\b)$ where $\phi=\phi_L:A\ra\hat{A}$ is
the symmetric isogeny associated with a rigidified  line  bundle
$L$ on $A$, $Z_{\phi,m}=(\phi,m\id_A)(A)\simeq A/\ker(\phi_m)$
where  $\phi_m=\phi|_{A_m}$,  $\b$ is obtained from $L^m$ by
descent (such $\b$ always exists if $m$ is odd,
since $\ker(\phi_m)$ is  isotropic
with respect to $e^{L^m}$).   Then
$K(Y,Z,T)\simeq Z$, $M(Y,Z,T)\simeq\b$, $Y\cap T=0$,
$Y\cap Z\simeq \ker(\phi)/\ker(\phi_n)$, and
$Z\cap T\simeq A_n/\ker(\phi_n)$. Hence,
if we take $W=\O_S$ we get
$$M[n]\simeq\und{\Hom}_{G_{Y,Z}\times_{\G_m}G_{Z,T}}(V\ot U,p_*\b).$$
In particular, when  $m=1$  we  have  $Z\simeq  A$,  $\b=L$,
$Z\cap T=0$, and $G_{Y,Z}=G(L)$. Thus, if we take $U=W=\O_S$
we obtain $M[n]\simeq\und{\Hom}_{G(L)}(V,p_*L)$.
      
Note  that if one of the pairwise intersections of $Y$, $Z$
and $T$ is trivial
then $K(Y,Z,T)$ is an abelian scheme over $S$. More precisely,
if say $Y\cap Z=0$ then $K(Y,Z,T)\simeq T$ and it is easy to
see  from  the  above  considerations  that
in this case we have an isomorphism of Heisenberg groups
\begin{equation}\label{GYZT}
G(M(Y,Z,T))\simeq G_{Y,Z}\times_{\G_m}G_{Z,T}\times_{\G_m}G_{T,Y}.
\end{equation}

\section{Weil representation on the derived category of an abelian scheme}
      
In this section the base scheme $S$ is always assumed to
noetherian, normal and connected.
Let $K$ denotes the  field  of  rational
functions on $S$.
      
\begin{lem}\label{extend}
Let $A$ and $A'$ be abelian schemes over $S$,
$A_K$ and $A'_K$ be their general  fibers  which  are  abelian
varieties over $K$. Then  the  restriction map
$$\Hom_S(A,A')\ra\Hom_K(A_K,A'_K):f\mapsto f|_K$$
is  an  isomorphism.  The  morphism $f$ is an isogeny if and
only if $f|_K$ is an isogeny.
\end{lem}
      
\Pf . The proof of the first assertion
is similar to the proof of the fact that an
abelian  scheme over a Dedekind scheme is a N\'eron model of
its generic fiber (see \cite{Neron} 1.2.8).
We have to check that any homomorphism $f_K:A_K\ra A'_{K}$
extends to a homomorphism $f:A\ra A'$. Let $\phi:A\ra A'$ be
the rational  map  defined  by  $f_K$.  Since $A$ is normal, by a
valuative criterion  of  properness  $\phi$  is  defined  in
codimension $\leq 1$. Let $V\sub A$ be a non-empty subscheme
over  which  $\phi$  is  defined.  Then since the projection
$p:A\ra S$ is flat and of finite presentation  it  is  open.
Thus,  $U=p(V)$  is  open,  and  $\phi_U:A_U\ra  A'_U$  is a
$U$-rational map in the terminology of \cite{Neron}.
By Weil's  theorem  (see  \cite{Neron}  4.4.1)  $\phi_U$  is
defined  everywhere, hence we get a homomorphism $f_U:A_U\ra
A'_U$ extending $f_K$. It remains to invoke Prop. I  2.9  of
\cite{FC} to finish the proof.
      
The part concerning  isogenies  can
be  proven  by  exactly the same argument as in \cite{Neron}
7.3 prop. 6: starting with an isogeny $f|_K:A_K\ra A'_K$ we
can find another isogeny $g|_K:A'_K\ra A_K$ such that the
composition $g|_Kf|_K$ is the multiplication by  an  integer
$l$ on $A_K$. By the first part we can  extend  $g_K$  to  a
homomorphism $g:A'\ra A$. This implies that   $gf=l_A$  ---  the
multiplication by $l$ morphism on $A$. It follows  that  the
restriction of $f$ to each fiber is an isogeny, hence $f$ is
an isogeny itself.
\ed
      
For an abelian scheme $A$  the group $\SL_2(A)$ is
defined  as  the subgroup of automorphisms of $\hat{A}\times
A$    preserving    the    line    bundle    $p_{14}^*\PP\ot
p_{23}^*\PP^{-1}$ on $(\hat{A}\times A)^2$. More explicitly,
if we write an automorphism of $\hat{A}\times A$ as a matrix
$g=\left( \matrix a_{11} & a_{12}\\ a_{21} & a_{22}
\endmatrix \right)$ where $a_{11}\in\Hom(\hat{A},\hat{A})$,
$a_{12}\in\Hom(A,\hat{A})$ etc., then
$g\in\SL_2(A)$  if  and  only  if  the  inverse automorphism
$g^{-1}$ is given by the matrix
$\left(\matrix \hat{a}_{22} & -\hat{a}_{12}\\
-\hat{a}_{21} & \hat{a}_{11} \endmatrix \right)$.
It  follows  from  Lemma \ref{extend} that when the base $S$ is
normal we have $\SL_2(A)\simeq\SL_2(A_K)$.
    
Now similarly to the classical picture one has  to  consider
the group of automorphisms of the Heisenberg extension $H(X)$
corresponding to $X=\hat{A}\times A$ with the structure
of enhanced symplectic biextension given by $\BB=p_{14}^*\PP$.
Namely, we define $\wt{\SL}_2(A)$  as  the  group  of
triples $g=(\bar{g},L^g,M^g)$ where
$\ov{g}=
\left( \matrix a_{11} & a_{12}\\ a_{21} & a_{22} \endmatrix \right)
\in\SL_2(A)$,
$L^g$ (resp. $M^g$) is a line bundle on $\hat{A}$ (resp. $A$)
rigidified along the zero section, such that
\begin{equation}
\phi_{L^g}=\hat{a}_{11}a_{21},\ \phi_{M^g}=\hat{a}_{22}a_{12},
\end{equation}
where  for  a  line  bundle  $L$ on an abelian scheme $B$ we
denote by $\phi_L:B\ra\hat{B}$ the symmetric homomorphism
corresponding to the symmetric biextension $\La(L)$.
The group law on $\wt{\SL}_2(A)$
is defined uniquely from the condition that the there is
an action of $\wt{\SL}_2(A)$
on the stack of Picard groupoids $\HH(A)$ such that
an element $g=(\bar{g},L^g,M^g)$ acts by the functor which is
identical on $\PPic$ and
sends the generator      $T_{(x,y)}$      (where
$(x,y)\in\hat{A}\times A$) to
$L_x\ot M_y\ot \PP_{(a_{12}y,a_{21}x)}\ot T_{\bar{g}(x,y)}$.
We refer to \cite{Weilrep} for explicit formulas for the
group law in $\wt{\SL}_2(A)$. It is easy  to  see  that  the
natural projection $\wt{\SL}_2(A)\ra\SL_2(A)$   is   a
homomorphism    with    the     kernel     isomorphic     to
$A(S)\times\hat{A}(S)$.
    
Consider the subgroup $\widehat{\SL}_2(A)\sub\wt{\SL}_2(A)$
consisting of triples with symmetric $L^g$ and $M^g$. Then we have
an isomorphism $\widehat{\SL}_2(A)\simeq\widehat{\SL}_2(A_K)$
since any symmetric  line  bundle  on  $A_K$  extends  to  a
symmetric line bundle on $A$ (see \cite{MB}, II.3.3).
      
Let $\Ga(A)=\Ga(A_K)$ be the image of the projection
$\widehat{\SL}_2(A_K)\ra\SL_2(A_K)$.
Then $\Ga(A)$ has finite index in $\SL_2(A_K)$
since it contains the subgroup $\Ga(A,2)=\Ga(A_K,2)\sub\SL_2(A_K)$
consisting of matrices with $a_{12}$ and $a_{21}$ divisible by 2.
    
In the case when $S$ is the spectrum of an algebraically closed field
it was shown in \cite{Weilrep} that there exist
intertwining functors between representations of Heisenberg
groupoid corresponding to the natural action
of  $\wt{\SL}_2(A)$    on    the    Heisenberg    groupoid
$H(\hat{A}\times A)$ which are analogous to the operators
of  Weil-Shale representation. We are going
to extend this construction to the case of  a  normal  base
scheme.
      
Recall (see \cite{Weilrep}, sect. 10) that there is a natural action of
$\wt{\SL}_2(A)$ on the set of enhanced lagrangian subvarieties
in    $X=\hat{A}\times    A$    such    that    a     triple
$g=(\bar{g},L^g,M^g)\in\wt{\SL}_2(A)$
maps $(\hat{A},\O_{\hat{A}})$ to  $\bar{g}(\hat{A})=
(a_{11},a_{21})(\hat{A})$
with the line bundle corresponding to $L^g\in\Pic(\hat{A})$.
Furthermore, there is a natural equivalence of categories
$$\ov{g}_*:\DD(\hat{A},\O_{\hat{A}})\ra
\DD(\ov{g}(\hat{A}),L^g):\FF\mapsto \ov{g}_*\FF$$
such that the standard $H(X)$-action
on $\DD(\hat{A},O_{\hat{A}})$ corresponds  to  the $g$-twisted
$H(X)$-action on $\DD(\ov{g}(\hat{A}),L^g)$. On the other hand,
if  $\hat{A}\cap  \ov{g}(\hat{A})=\ker(a_{21})$ is finite
(hence, flat) over $S$ and there   exists   a
Schr\"odinger  representation $V$ for  the corresponding
Heisenberg extension
$G_g:=G_{\hat{A},\ov{g}(\hat{A})}$ of $\ker(a_{21})$ then
the   construction   of   the   previous  section  gives
another equivalence
$$\DD(\ov{g}(\hat{A}),L^g))\ra\DD(\hat{A},\O_{\hat{A}})$$
compatible with the standard $H(X)$-actions. Composing it
with the previous equivalence we get an equivalence
$\rho\wt{\ra}\rho^g$ where $\rho$ is the representation
of $H(X)$ on $\DD(\hat{A},\O_{\hat{A}})\simeq\DD(A)$
given by (\ref{act}),
$\rho^g=\rho\circ  g$  is the same representation twisted by
$g$. Using (\ref{kernel})
it  is  easy to compute that the kernel on $A\times_S A$
corresponding to this equivalence has form
\begin{equation}\label{kernelnew}
\KK(g,V)=(p_2-a_{22}p_1)^*E\ot (a_{12}\times\id)^*\PP^{-1}\ot
(-p_1)^*M^g
\end{equation}
where
$$E=\und{\Hom}_{G_g^{-1}}(V^*,a_{21*}((L^g)^{-1})).$$
Note that here $G_g^{-1}$ is the restriction of the
Mumford's extension $G((L^g)^{-1})\ra\ker(\hat{a}_{11}a_{21})$
to $\ker(a_{21})$.
 
Hence,
if $\ker(a_{21})$ is finite over $S$  we get a functor from the
gerb of Schr\"odinger representations for $G_g$ to the stack of
$H(X)$-equivalences   $\Isom_{H(X)}(\rho,\rho^g)$,
More precisely, the category of intertwining operators between
$\rho$ and $\rho^g$ is defined in terms of kernels in
$\D^b(A\times_S A)$ (see \cite{Weilrep}) and the glueing property
is satisfied because the kernels corresponding to
equivalences are actually vector bundles (perhaps, shifted).
Indeed, the latter property is local  with  respect  to  the
{\it fppf} topology on the base $S$ and locally a Schr\"odinger
representation for $G_g$ exists and give rise to
the  kernel (\ref{kernelnew}) in  $\Isom_{H(X)}(\rho,\rho^g)$  which is a
vector bundle up to shift. Now any other object of
$\Isom_{H(X)}(\rho,\rho^g)$ is obtained from a given one
by tensoring with a line bundle on $S$ and a shift.
      
Thus, when $\ker(a_{21})$ is finite the
obstacle for the existence of a global equivalence
between   $\rho$ and   $\rho^g$  is  given  by  the  class
$e(G_g)\in\Br(S)$.
Let $U\sub\SL_2(A)$ be the subset of matrices such that
$a_{21}$ is an isogeny. It turns out that similarly to the
case of real groups one can deal with $U$ instead of
the entire group when defining representation of $\SL_2(A)$.
This observation can be formalized as follows.
Let us call a subset $B$ of a group
$G$  {\it  big}  if  for  any  triple of elements $g_1, g_2,
g_3\in G$ the intersection $B^{-1}\cap Bg_1\cap Bg_2\cap Bg_3$
is non-empty. This condition first appeared
in \cite{Weil} IV. 42, while the term is due to D.~Kazhdan.
The reason for introducing this notion is the following lemma.
   
\begin{lem} Let $B\sub G$ be a big subset.
Then $G$ is isomorphic  to
the abstract group generated by elements $[b]$ for $b\in B$
modulo  the relations $[b_1][b_2]=[b]$ when $b,b_1,b_2\in B$
and $b=b_1b_2$.
If $c:G\times G\ra C$
is a 2-cocycle (where $C$ is an abelian group with the trivial
$G$-action) such that $c(b_1,b_2)=0$ whenever
$b_1, b_2, b_1b_2\in B$ then $c$ is a coboundary.
\end{lem}
   
\Pf . For the proof of the first statement we refer to
\cite{Weil}, IV, 42, Lem. 6.
Let $c:G\times G\ra C$ be a 2-cocycle,
$H$ be the corresponding central extension of $G$ by $C$.
Consider the group $\wt{H}$ generated by the central
subgroup $C$ and generators $[b]$ for $b\in B$
subject to relations $[b_1][b_2]=c(b_1,b_2)[b_1b_2]$,
where $b_1, b_2, b_1b_2\in B$. Then
$\wt{H}/C\simeq G$, hence the natural homomorphism
$\wt{H}\ra H$ is an isomorphism. If $c(b_1,b_2)=0$
whenever $b_1, b_2, b_1b_2\in B$, then the
extension $\wt{H}\simeq H\ra G$ splits, hence
$c$ is a coboundary.
\ed
   
At this point we need to recall
some results from \cite{Weilrep}, sect. 9 concerning
the group $\SL_2(A)$.
Since $\SL_2(A)=\SL_2(A_K)$ we can work with abelian varieties
over a field. First note that this group
can be considered as a group of $\Z$-points of an
group scheme over $\Z$. It turns out that the corresponding
algebraic group $\SL_{2,A,\Q}$ over $\Q$ is very close to be
semi-simple. Namely, if we fix a polarization on $A$ then the
latter  group  is  completely  determined  by  the   algebra
$\End(A)\ot\Q$ and the Rosati involution on it. Decomposing
(up to isogeny) $A$ into a product $A_1^{n_1}\times\ldots
A_l^{n_l}$ where $A_i$ are different simple abelian varieties
and choosing a polarization compatible with this decomposition
it is easy to see that
$$\SL_{2,A,\Q}\simeq\prod_i R_{K_{i,0}/\Q}\U^*_{2n_i,F_i}$$
where $F_i=\End(A_i)\ot\Q$, $K_i$ is the center of $F_i$,
$K_{i,0}\sub K_i$ is the subfield of elements stable under
the  Rosati  involution,  $\U^*_{2n_i,F_i}$  is  the group of
$F_i$-automorphisms of $F_i^{2n_i}$
preserving the standard skew-hermitian form,
$R_{K_{i,0}/\Q}$ denotes the restriction of scalars from
$K_{i,0}$ to $\Q$. Thus, the only case
when the group $\U^*_{2n_i,F_i}$ is not semi-simple is when
the Rosati involution on $F_i$ is of the second kind, i.~e.
$K_{i,0}\neq K_i$. In the latter case, $\U^*_{2n_i,F_i}$
is a product of the semi-simple subgroup $\SU^*_{2n_i,F_i}$
(defined using the determinant with values in $K_i$)
and the central subgroup
$K^1_i=\{ x\in K_i\ |\ N_{K_i/K_{i,0}}=1\}$ consisting
of diagonal matrices. Furthermore, the intersections of these
two subgroup is finite. It follows that the group
$\SL_{2,A,\Q}$ always has an almost direct decomposition into
a product of the semi-simple subgroup
$H=\prod_i R_{K_{i,0}/\Q}\SU^*_{2n_i,F_i}$
and a central subgroup $Z$ consisting of diagonal matrices.
Now we can prove the following result.
   
\begin{lem}\label{big}
Let $\Ga\sub\SL_2(A)$ be a subgroup of finite
index. Then the subset $\Ga\cap U$ is big.
\end{lem}
      
\Pf . By definition $U$ is an intersection of a Zariski open
subset $\und{U}$ in the irreducible algebraic group
$\SL_{2,A_K,\Q}$  with  $\Ga$. If $\Ga$
were Zariski dense in $\SL_{2,A_K,\Q}$ the  proof  would  be
finished. This is not always true, however, we claim that
$\Ga\cap H$ is dense in $H$ where $H$ is the subgroup of
$\SL_{2,A_K,\Q}$ introduced above. Indeed, this follows from
the fact that $H$ is semi-simple and the corresponding
real groups have no compact factors (see \cite{Weilrep} 9.4).
On the other hand, the set $\und{U}$ is $Z$-invariant (since
$Z$ consists of diagonal matrices).
It follows that for any $g\in\Ga$
the intersection $\und{U}g\cap H$ is a non-empty  Zariski
open subset of $H$ (as a preimage of a non-empty Zariski set
under the isogeny $H\ra\SL_{2,A_K,\Q}/Z$). Therefore, for
any triple of elements $g_1,g_2,g_3\in\Ga$ the intersection
$\und{U}\cap\und{U}g_1\cap\und{U}g_2\cap\und{U}g_3\cap H$
is a non-empty  open  subset  in  $H$,  hence  it contains an
element of $\Ga\cap H$. As $U=U^{-1}$ this shows that $U$ is
big.
\ed

\begin{prop}\label{obshom}
There exists a homomorphism
$\de_A:\Ga(A)\ra\Br(S)$  such that for any $g\in\widehat{\SL}_2(A)$
lying over $\ov{g}\in\Ga(A)$ there exists a global object in
$\Isom_{H(X)}(\rho,\rho^g)$ if and only if $\de_A(\ov{g})=0$.
If the $a_{21}$-entry of $\ov{g}$ is an isogeny then
$\de_A(\ov{g})=e(G_g)$.
\end{prop}
      
\Pf . For any $g\in\widehat{\SL}_2(A)$ let us denote 
by $\Isom_{H(X)}^0(\rho,\rho^g)\sub\Isom_{H(X)}(\rho,\rho^g)$
the full subcategory of kernels that belong to the core of
the standard $t$-structure on $\D^b(A\times_S A)$. We claim
that $\Isom_{H(X)}^0(\rho,\rho^g)$ is a gerb.
Indeed,  we  know  this when $\ov{g}\in\Ga(A)\cap U$. Now by
Lemma \ref{big} any element
in $\widehat{\SL}_2(A)$ can be written as a product
$gg'$ where $g$ and $g'$ lie over $\Ga(A)\cap U$.
Locally over $S$
there exist Schr\"odinger representations  for  $G_g$  and
$G_{g'}$,  hence  (\ref{kernelnew})  defines  the  corresponding
kernels $\KK(g)\in\Isom_{H(X)}(\rho,\rho^g)$ and
$\KK(g')\in\Isom_{H(X)}(\rho,\rho^{g'})$. Now their composition
$\KK(gg')=p_{13*}(p_{12}^*\KK(g')\ot p_{23}^*\KK(g))$ is an element of
$\Isom_{H(X)}(\rho,\rho^{gg'})$ and  Prop.  1.7  of  \cite{Muk2}
implies that $\KK(gg')$ has only one non-zero sheaf cohomology
which is flat over $S$
(since this is so in the  case  of  an algebraically closed field
considered in \cite{Weilrep}). It follows that any object of
$\Isom_{H(X)}(\rho,\rho^{gg'})$ over any open subset of $U$
is a pure $S$-flat sheaf (perhaps shifted), hence,  the  gluing  axiom  is
satisfied.
      
Now we can  define  $\de_A(\ov{g})\in H^2(S,\G_m)$ for  an
element $\ov{g}\in\Ga(A)$ as the class of the gerb
$\Isom_{H(X)}^0(\rho,\rho^g)$ where $g\in\widehat{\SL}_2(A)$ is
any element lying over $\ov{g}$.
When $\ov{g}\in U$ this class is equal to $e(G_g)\in\Br(S)$.
Clearly,  $\de_A$ is a homomorphism $\Ga(A)\ra H^2(S,\G_m)$.
Since $\Ga(A)$ is generated by $\Ga(A)\cap U$ we have
$\de_A(\ov{g})\in\Br(S)$ for any $\ov{g}\in\Ga(A)$.
\ed

\begin{lem}\label{neutrcomp}
Let   $g=(\ov{g},L^g,M^g)$   be   an   element   of
$\widehat{\SL}_2(A)$  such that $\ker({a_{21}})$ is flat over
$S$, its neutral component $\ker(a_{21})^0$ is an abelian
subscheme of $\hat{A}$, and $L^g|_{\ker(a_{21})^0}$ is trivial.
Then $\de_A(\ov{g})=0$ if and only if there exists a  Schr\"odinger
representation      for     the     Heisenberg     extension
$G_g$ of
$\pi_0(\ker(a_{21}))=\ker(a_{21})/(\ker(a_{21}))^0$ induced by
$G(L^g)|_{\ker(a_{21})}$.
\end{lem}
      
\Pf . Let $V$ be a Schr\"odinger representation for
$G_g$. Then we can define $\KK(g,V)$ by the formula (\ref{kernelnew})
where
$E$ is defined as follows. First we descend $L^g$ to a
line bundle $\ov{L}$ on $\hat{A}/\ker(a_{21})^0$, then
we set
$$E=\und{\Hom}_{G_g^{-1}}(V^*,\ov{a}_{21*}(\ov{L}^{-1}))$$
where $\ov{a}_{21}:\hat{A}/\ker(a_{21})^0\ra A$ is the
finite map induced by $a_{21}$. Note that when  $S$  is  the
spectrum of an algebraically closed field this kernel
$\KK(g,V)$ coincides with the one defined in \cite{Weilrep},
12.3. By definition $\KK(g,V)$
is the direct image of a bundle on an abelian subscheme
$$\supp  \KK(g,V)=\{(x_1,x_2)\in   A^2\   |   x_2-a_{22}x_1\in
a_{21}(\hat{A})\}$$
(note  that  $a_{21}(\hat{A})\sub A$ is an abelian subscheme
since $\ker(a_{21})$ is flat). Applying this to $g^{-1}$  we
get
$$\supp      \KK(g^{-1},V^*)=\{(x_1,x_2)\in      A^2\     |\
x_2-\hat{a}_{11}x_1\in
\hat{a}_{21}(\hat{A})\}.$$
Hence, the sheaf $p_{12}^*\KK(g^{-1},V^*)\ot p_{23}^*\KK(g,V)$ on
$A^3$ is supported on the abelian subscheme
$$X(g)=\{(x_1,x_2,x_3)\in     A^3\    |\    x_3-a_{22}x_2\in
a_{21}(\hat{A}),
x_2-\hat{a}_{11}x_1\in\hat{a}_{21}(\hat{A})\}.$$
Note that for $(x_1,x_2,x_3)\in X(g)$ we have
$$x_3-a_{22}\hat{a}_{11}x_1\in a_{21}(\hat{A})+
a_{22}\hat{a}_{21}(\hat{A})=a_{21}(\hat{A})$$
since $a_{22}\hat{a}_{21}=a_{21}\hat{a}_{22}$. Now
$a_{22}\hat{a}_{11}x_1\equiv x_1\mod(a_{21}(\hat{A}))$,
hence $x_3\equiv x_1\mod(a_{21}(\hat{A}))$.
Thus, we have an isomorphism
\begin{align*}
&X(g)\ra\{(x,x')\ |\  x-x'\in a_{21}(\hat{A})\}
\times\hat{a}_{21}(\hat{A}):\\
&(x_1,x_2,x_3)\mapsto ((x_1,x_3), x_2-\hat{a}_{11}x_1).
\end{align*}
It follows that the restriction of
the projection $p_{13}$ to $X(g)$ is flat and surjective.
Therefore, applying
Prop. 1.7 of \cite{Muk2} we conclude as in the proof of Theorem
\ref{eq} that $\ov{R}(V)$ is an equivalence.
Thus, the map $V\mapsto \KK(g,V)$ gives a functor from
$\Sch_{G_g}$ to $\Isom_{H(X)}(\rho,\rho^g)$.
\ed

\begin{prop}\label{compat}
Let $A$ and $B$ be  abelian  schemes  over  the
same base $S$, then there is a natural embedding
$i_A:\Ga(A)\ra\Ga(A\times B)$ such that
$$\de_A=\de_{A\times B}\circ i_A.$$
\end{prop}
      
\Pf . It is sufficient to check this identity on elements
of $\Ga(A)\cap U$, in which case  this  follows  immediately
from Lemma \ref{neutrcomp}.
\ed
      
Consider the subgroup of finite index
$\Ga_0(A)\sub\Ga(A)$ consisting of matrices
$\left( \matrix a_{11} & a_{12}\\ a_{21} & a_{22} \endmatrix \right)$
in $\Ga(A)$ for which $a_{21}$ is divisible by 2.
Note that $\Ga_0(A)$ contains $\Ga(A,2)$.
      
\begin{lem}\label{oddisog}
Let $A$ be an abelian variety over a field such
that there exists a symmetric line bundle $L$ on $A$
which induces an isogeny $f:A\ra\hat{A}$ of odd degree.
Then the subgroup of $\Ga(A)$ generated by its elements
$\left( \matrix a_{11} & a_{12}\\ a_{21} & a_{22} \endmatrix \right)$,
for  which  $a_{21}$  is  an isogeny of odd degree, contains
$\Ga_0(A)$.
\end{lem}
      
\Pf. Recall that $U\sub\SL_2(A)$ is the subset defined by the
condition that $a_{21}$ is an isogeny. Consider the matrix
$\ga_f=\left(   \matrix   \id  &  0  \\  f  &  \id  \endmatrix
\right)\in\Ga(A)$.  Let  $U_1=U\cap U\ga_f^{-1}$.
Then the argument similar to that of Lemma  \ref{big}  shows
that $\Ga_0(A)\cap U_1$ is a big subset in $\Ga_0(A)$, in particular,
$\Ga_0(A)$  is  generated  by
$\Ga_0(A)\cap U_1$. Now let $\ga\in\Ga_0(A)\cap U_1$, then
$\ga \ga_f\in U$ and its $a_{21}$-entry is an isogeny of odd
degree which implies the statement.
\ed

\begin{prop}\label{triv1}
The restriction of the homomorphism
$\de_A:\Ga(A)\ra\Br(S)$ to $\Ga_0(A)$ is trivial.
\end{prop}
      
\Pf . Recall that if $a_{21}$ is an isogeny then
$\de_A(g)$ is defined as an obstacle for
the existence of a Schr\"odinger representation of the
(symmetric) Heisenberg extension $G_g$ of $\ker(a_{21})$ attached to $g$,
where $g$ projects to a matrix
$\left( \matrix a_{11} & a_{12}\\ a_{21} & a_{22} \endmatrix
\right)\in\SL_2(A)$. Hence, by Theorem \ref{odd} $\de_A(g)=0$
if $a_{21}$ is an isogeny of odd degree. In particular,
if $A$ is principally polarized then Lemma \ref{oddisog}
implies that the restriction of $\de_A$ to $\Ga_0(A)$
is trivial.
      
By Zarhin's trick (see \cite{Zarh}) for any abelian scheme $A$ over $S$
there exists an abelian scheme $B$ over $S$ such that
$A\times   B$   admits  a  principal  polarization.  Now  by
Proposition \ref{compat} we have
$\de_A=\de_{A\times B}\circ i_A$. Therefore,
the restriction of $\de_A$ to $\Ga_0(A)$ is trivial.
\ed
      
\begin{rem} It is easy to see that the kernel of $\de_A$ is in
general bigger than $\Ga_0(A)$. Namely, it contains also matrices
for which $a_{21}$ (or $a_{12}$) is an isogeny of odd degree,
those for which $a_{12}$ is divisible by 2, and those for which
$a_{11},a_{22}\in\Z$. Sometimes,
these elements together with $\Ga_0(A)$ generate the entire
group $\Ga(A)$,
however, it is not clear whether $\de_A$ is always trivial.
In the section  \ref{real} we will prove it in some special cases.
\end{rem}
      
Let $\widehat{\Ga}_0(A)$ be the preimage of $\Ga_0(A)$ in
$\widehat{\SL}_2(A)$. In other words, this is the subgroup
of elements $g\in\widehat{\SL}_2(A)$ such that for the
corresponding matrix in $\SL_2(A)$ the $a_{21}$-entry
is divisible by 2.
We say that there is a faithful action of a group
$G$ on a category $\CC$ if there is an embedding of $G$ into
a  group  of  autoequivalences  of  $\CC$  (considered up to
isomorphism).
      
\begin{thm} For any abelian scheme $A$ over a normal connected
noetherian base $S$ there is a faithful action of a central extension of
the group $\widehat{\Ga}_0(A)$ by $\Z\times\Pic(S)$ on $\D^b(A)$.
\end{thm}
   
\Pf . According to Proposition \ref{triv1}
for every $g\in\widehat{\Ga}_0(A)$ there
exists a global object in $\Isom_{H(X)}(\rho,\rho^g)$.
It is defined uniquely up to a shift and tensoring with a
line bundle on $S$. Hence, the required action of a central
extension. The fact that this action is faithful is
clear in the case when the base is a field: for example,
one can use explicit formulas for these functors from
Lemma \ref{neutrcomp}. Since  the  action  of  $\Pic(S)$  on
$\D^b(A)$ is obviously faithful the general case follows.
\ed

\begin{cor} With the assumptions of the above theorem, there
is a faithful action of a central extension of $\Ga(A,2)$
by $\Z\times\Pic(S)$ on $\D^b(A)$.
\end{cor}
   
\Pf . This action is obtained from the canonical homomorphism
$\Ga(A,2)\ra\widehat{\SL}_2(A)$ splitting
the projection $\widehat{\SL}_2(A)\ra\SL_2(A)$. Namely,
under this splitting the matrix
$\ov{g}=\left( \matrix a_{11} & a_{12}\\ a_{21} & a_{22} \endmatrix
\right)\in\Ga(A,2)$ maps to the element
$(\ov{g},(\hat{a}_{11}a_{21}/2,\id)^*\PP,
(\id,\hat{a}_{22}a_{12}/2)^*\PP)$ of $\widehat{\SL}_2(A)$,
where $\PP$ is the Poincar\'e line bundle on $A\times\hat{A}$.
\ed
   
  
\section{The induced action on a Chow motive}
                                        
In this section we will construct a projection action of
the algebraic group $\SL_{2,A,\Q}$ on the relative Chow motive
of an abelian scheme $\pi:A\ra S$ with rational coefficients.
Let us denote by $\Cor(A)$ the Chow group
$\CH^*(A\times_S A)\ot\Q$
considered as  a
$\Q$-algebra with multiplication given by the composition of
correspondences:
\begin{equation}\label{compcorr}
\beta\circ\a=p_{13*}(p_{12}^*(\a)\cdot p_{23}^*(\beta))
\end{equation}
where $\a,\beta\in\CH^*(A\times_S A)\ot \Q$,
$p_{ij}$ are the projections from $A^3$ to $A^2$.
The       unit      of      this    algebra     is
$[\De]\in\CH^g(A\times_S A)$ where $\De\sub A\times_S A$ is the
relative diagonal, $g=\dim A$.
   
Using the Riemann-Roch theorem
it is easy to see that the multiplication (\ref{compcorr})
is compatible with
the composition law on $K^0(A\times_S A)$ arising from the
interpretation  of  $\D^b(A\times_S A)$  as  the  category of
functors from $\D^b(A)$ to itself considered above,
via the map
\begin{equation}\label{ch}
K^0(A\times_S A)\ot \Q\wt{\ra}\CH^*(A\times_S A)\ot\Q:
x\mapsto \ch(x)\cdot\pi^*\Td(e^*T_{A/S})
\end{equation}
where $\ch$ is the Chern character.
Let us consider
the embedding of algebras
$$\CH(S)_{\Q}\ra\Cor(A):x\mapsto \pi^*(x)\cdot [\De]$$
where $\CH(S)_{\Q}=\CH(S)\ot\Q$ is equipped with the
usual multiplication. In particular,
we have an embedding of groups of invertible elements
$\CH(S)_{\Q}^*\sub\Cor(A)^*$.
Applying the map (\ref{ch}) to the kernels giving
projective action of $\Ga(A,2)$ on $\D^b(A)$ we obtain
a homomorphism
$$\wt{\phi}:\Ga(A,2)\ra(\Cor(A))^*/\pm\Pic(S),$$
where $\Pic(S)$ is embedded into $\CH(S)_{\Q}^*$
by Chern character (multiplication by $\pm$
arises from shifts in derived category).
Our aim is to approximate this homomorphism by
a morphism of algebraic groups over $\Q$.
More precisely, we have to replace $\wt{\phi}$ by the induced
homomorphism
$$\phi:\Ga(A,2)\ra(\Cor(A))^*/\CH(S)_{\Q}^*.$$
Now we claim that one can replace here source and target
by some algebraic groups over $\Q$ such that $\phi$ will be
induced by an algebraic homomorphism.
Naturally, the source should be replaced by
$\SL_{2,A,\Q}$ (see the previous section). To
approximate the target we have to replace algebras $\Cor(A)$ and
$\CH(S)_{\Q}$ by their finite-dimensional subalgebras.
   
\begin{thm}\label{actmotmain}
There exists a finite-dimensional $\Q$-subalgebra
$D\sub\Cor(A)$ and a morphism of algebraic
$\Q$-groups $\rho:\SL_{2,A,\Q}\ra D^*/(D\cap\CH(S)_{\Q})^*$
inducing $\phi$ on $\Ga(A,2)$.
\end{thm}
         
\Pf . For a pair of abelian schemes $A$ and $B$ over $S$
let us consider the map
$$\ga_{A,B}:\Hom(A,B)\ra\CH(A\times_S B)$$
that sends an $S$-morphism $f:A\ra B$ to
the class $[\Ga_f]$ of the (relative) graph of $f$.
One can extend naturally $\ga$ to a
map
$$\ga:\Hom(A)\ot\Q\ra\CH(A\times_S B)\ot\Q$$
by sending
$f/n$ to $\ga([n]_A)^{-1}\circ\ga(f)$, where
$f\in\End(A)$, $n\neq0$ --- here we take the inverse
to $\ga([n]_A)$ in the algebra $\Cor(A)$ and
use its natural action on $\CH(A\times_S B)\ot\Q$.
It is easy to check that $\ga$ is a polynomial map
(see \cite{Weilrep}, Lemma 13.3, for the case $A=B$).
  
Now let
$\ov{g}=\left( \matrix a_{11} & a_{12}\\ a_{21} & a_{22} \endmatrix
\right)$  be an element of $\Ga(A,2)\cap U$ (recall that $U$
is defined by the condition that $a_{21}$ is an isogeny).
From the formula (\ref{kernelnew}) we
get the following expression for $\phi(\ov{g})$:
\begin{equation}\label{chker}
\phi(\ov{g})=(p_2-a_{22}p_1)^*(a_{21*}\ch(L^g)^{-1})\cdot
(a_{12}\times\id)^*(\ch(\PP))\cdot
p_1^*(\ch(M^g))\mod\CH(S)_{\Q}^*.
\end{equation}
Note that $\ch(L^g)$ (resp. $\ch(M^g)$) is a polynomial function of
$\hat{a}_{11}a_{21}$   (resp.   $\hat{a}_{22}a_{12}$).
Also  the  functors  $f^*$  and  $f_*$  can  be expressed as
compositions with the correspondence given by the graph of $f$.
It follows from the above remarks that the right hand side of (\ref{chker})
is obtained by evaluating at $\ov{g}$ of a polynomial map
$\psi:\End(\hat{A}\times A)\ot\Q\ra\Cor(A)$.
In particular, the image   of   $\psi$  belongs  to  a
finite-dimensional $\Q$-subspace of $\Cor(A)$.
  
Let $\und{U}\sub\SL_{2,A,\Q}$ be the Zariski open
subset defined by $\deg(a_{21})\neq0$. Note that
$\und{U}$  is  stable under the inversion morphism $g\mapsto
g^{-1}$. Let us also denote
$\und{U}^{(2)}=
\mu^{-1}(\und{U})\cap(\und{U}\times\und{U})\sub\und{U}\times
\und{U}$ where $\mu$ is the group law.
Consider two polynomial maps
\begin{align*}
&a_1:\und{U}\ra\Cor(A):u\mapsto a_1(u)=\psi(u^{-1})\circ\psi(u),\\
&a_2:\und{U}^{(2)}\ra\Cor(A):
(u_1,u_2)\mapsto a_2(u_1,u_2)=\psi(u_1u_2)\circ\psi(u_2^{-1})\circ
\psi(u_1^{-1}).
\end{align*}
It is easy to see that the images of both  maps  belong  to  the
subalgebra $\CH(S)_{\Q}\sub\Cor(A)$. This can be done either
by direct computation using (\ref{chker}) or using the density
of $\Ga$ in $H\sub\SL_{2,A,\Q}$ (see the  previous  section).
Also an easy direct computation shows that
$a_1(u)$ is invertible in the algebra $\CH(S)_{\Q}$ for all
$u\in\und{U}$. This immediately implies that $a_2(u_1,u_2)$
is invertible for any $(u_1,u_2)\in\und{U}^{(2)}$: indeed,
$a_2(u_1,u_2)$ is a divisor of
$a_1(u_1)a_1(u_2)a_1(u_2^{-1}u_1^{-1})$.   Note   that   the
components of images of $a_1$ and $a_2$ span finite-dimensional subspaces
in $\CH(S)^i_{\Q}$ for any $i$. It follows that there exists a
finite dimensional subalgebra $D_S\sub\CH(S)_{\Q}$ such that
the images of $a_1$ and $a_2$ belong to $D_S^*$. Now we have
\begin{equation}\label{projmot}
\psi(u_1)\circ\psi(u_2)=a_2(u_1,u_2)^{-1}a_1(u_1)a_1(u_2)
\psi(u_1u_2)
\end{equation}
for $(u_1,u_2)\in\und{U}^{(2)}$. Let $D\sub\Cor(A)$ be the
$D_S$-submodule generated by $\psi(u)$ with $u\in\und{U}$.
Then $D$ is finite-dimensional as a $\Q$-vector space
and (\ref{projmot}) shows that
$\psi(u_1)\circ\psi(u_2)\in D$ for any
$(u_1,u_2)\in\und{U}^{(2)}$. Since $\und{U}^{(2)}$ is dense
in $\und{U}\times\und{U}$ it follows that $D$ is a subalgebra.
Now (\ref{projmot}) implies that $\psi$ uniquely extends to
a homomorphism $\SL_{2,A,\Q}\ra D^*/D_S^*$.
\ed
  
\section{Splittings  of the extension
$\wt{SL}_2(A)\ra SL_2(A)$}

Let $A/S$ be an abelian scheme with a principal
polarization $\phi:A\ra\hat{A}$. Then we have the
Rosati involution
$$\e_{\phi}:\End(A)\ra\End(A):
f\mapsto \phi^{-1}\circ \hat{f}\circ \phi.$$
The group $\SL_2(A)$ is completely determined by
the algebra $\End(A)$ with involution $\e_{\phi}$.
The definition of the group $\wt{\SL}_2(A)$
requires in addition the knowledge of the extension
\begin{equation}\label{extA}
0\ra \hat{A}(S)\ra\Pic(A)\ra\Hom^{\sym}(A,\hat{A})
\end{equation}
together with the action of the multiplicative monoid of $\End(A)$ on it.
Thus, splittings of the homomorphism
$\wt{\SL}_2(A)\ra\SL_2(A)$ should be related to
splittings of (\ref{extA}).
More precisely, it's natural to consider 
splittings compatible with the $\End(A)$-action. We'll show that such
splittings of (\ref{extA}) correspond to simultaneous splittings
of homomorphisms $\wt{\SL}_2(A^n)\ra\SL_2(A^n)$ for all $n$,
where $A^n/S$ is the $n$-th relative cartesian power of $A/S$.

More generally, we start with arbitrary subring $R\sub\End(A)$
stable under the Rosati involution. Let us denote by
$\e:R\ra R$ the restriction of the Rosati involution to $R$,
let $R^+\sub R$ be the subring of elements stable under $\e$.
Then for any $n\ge 1$ we can consider
the subgroup  $\SL_2(A^n,R)\sub\SL_2(A^n)$  consisting  of
$2n\times 2n$ matrices
with all entries belonging to $R$ (we identify $\hat{A}$ with
$A$ via $\phi$). Let $\wt{\SL}_2(A^n,R)\sub\wt{\SL}_2(A^n)$
be the preimage of $\SL_2(A^n,R)$. We are interested in splittings
of the natural homomorphisms
\begin{equation}\label{homSL}
\wt{\SL}_2(A^n,R)\ra\SL_2(A^n,R)
\end{equation}
It turns out that the following structure on $A$ is relevant for
this.

\begin{defi} A $\Sigma_{R,\e}$-structure for $\phi$ is a
homomorphism $R^+\ra\Pic(A):r_0\mapsto L(r_0)$
such that
\begin{equation}\label{CM1}
\phi_{L(r_0)}=\phi\circ [r_0]_A
\end{equation}
for any $r_0\in R^+$ and
\begin{equation}\label{CM2}
[r]^*L(r_0)\simeq L(\e(r)r_0r)
\end{equation}
for any $r\in R$, $r_0\in R^+$.
\end{defi}

Note that (\ref{CM2}) for $r=-1$ implies that all line bundles
$L(r_0)$ are symmetric.
In \cite{thetaid} we studied the question of existence of
$\Sigma_{R,\e}$-structure for an
abelian variety.
For example, in the case of a complex elliptic curve $E$ with  complex
multiplication such a structure for $R=\End(E)$ exists
if and only if $R$ is unramified at $2$. Another example
is the case when $R$ is a ring of integers in a
totally real number field unramified at 2 and $\e=\id$
(see next section).

\begin{thm}\label{Sig}
A $\Sigma_{R,\e}$-structure on $A$
induces canonical splittings of the homomorphisms
(\ref{homSL}) for all $n$.
\end{thm}

\Pf . It is easy to see that a $\Sigma_{R,\e}$-structure
on $A$ induces a similar structure on $A^n$ with $R$
replaced by the matrix algebra $\Mat_n(R)$ and
$\e$ replaced by the corresponding
involution $(a_{ij})\mapsto (\e(a_{ji}))$ of $\Mat_n(R)$.
Hence, it is sufficient to consider the case $n=1$.
In this case we define the splitting
$$\SL_{2}(A,R)\ra\widehat{\SL}_2(A,R):
g=\left( \matrix a_{11} & a_{12}\\ a_{21} & a_{22}  \endmatrix
\right)\mapsto (g,
(\phi^{-1})^*L(\e(a_{11})a_{21}),
L(\e(a_{22})a_{12}))$$
\ed

Now we are going to prove that conversely the existence
of splitting of (\ref{homSL}) for $n=2$ implies the
existence of $\Sigma_{R,\e}$-structure. We use the following
observation. For any abelian scheme $A$ there is a natural embedding of
the semi-direct  product
$\Aut(A)\ltimes\Hom^{\sym}(A,\hat{A})$ into $\SL_2(A)$ as the
subgroup of matrices with $a_{21}=0$. Thus, a splitting of
the homomorphism $\wt{\SL}_2(A)\ra\SL_2(A)$ restricts to
a splitting of the homomorphism of $\Aut(A)$-modules
$\Pic(A)\ra\Hom^{\sym}(A,\hat{A})$.   In    the    situation
considered above we have similar subgroups in
$\SL_n(A^n,R)$. Note that  the subgroup of
$\Hom^{\sym}(A^n,\hat{A}^n)$ consisting of matrices with
entries in $R$ can be identified with the group of
hermitian matrices $\Mat^{\herm}_n(R)$ and the natural right
action of $\GL_n(R)$ on it induced by its action on $A^n$
is given by the formula
$$B\mapsto \ov{C}BC$$
where $C=(c_{ij})\in\GL_n(R)$, $B\in\Mat^{\herm}_n(R)$,
$\ov{C}=(\e(c_{ji}))$.
Thus,  a  splitting  of  the  homomorphism
(\ref{homSL}) induces a splitting of the homomorphism of
$\GL_n(R)$-modules
$$\Pic(A^n,R)\ra\Mat^{\herm}_n(R)$$
where $\Pic(A^n,R)\sub\Pic(A^n)$ is the subgroup
of line bundles $L$ such that
$\phi_L\in\Hom^{\sym}(A^n,\hat{A}^n)$ has entries in $R$.

Now we claim that such a splitting for $n=2$ leads to
a $\Sigma_{R,\e}$-structure on $A$.
 
\begin{thm} Any splitting of the homomorphism
of $\GL_2(R)$-modules
$\Pic(A^2,R)\ra\Mat^{\herm}_2(R)$ is
induced by a unique $\Sigma_{R,\e}$-structure.
\end{thm}

\Pf . Let
$$s:\Mat^{\herm}_2(R)\ra\Pic(A^2,R)$$
be such a splitting. Then for $r_0\in R^+$ one has
\begin{equation}\label{s1}
s\left( \matrix r_0 & 0 \\ 0 & 0  \endmatrix
\right)=p_1^*L(r_0)\otimes p_2^*\eta(r_0)
\end{equation}
for some line bundle $L(r_0)$ and $\eta(r_0)$ on $A$ such that
$\phi_{L(r_0)}=\phi\circ [r_0]$, $\eta(r_0)\in\Pic^0(A)$.
The compatibility of $s$ with the action of $\GL_2(R)$
means that
$$[C]^*s(B)=s(\ov{C}BC)$$
where $B\in\Mat^{\herm}_2(R)$, $C=(c_{ij})\in\GL_2(R)$,
$\ov{C}=(\e(c_{ji}))$. Applying this to
$C=\left( \matrix 0 & 1 \\ 1 & 0  \endmatrix
\right)$ we deduce from (\ref{s1}) the equality
\begin{equation}\label{s2}
s\left( \matrix 0 & 0 \\ 0 & r_0  \endmatrix
\right)=p_2^*L(r_0)\ot p_1^*\eta(r_0).
\end{equation}
Also using the identity
$$\left( \matrix 1 & 0 \\ \e(r) & 1  \endmatrix
\right)\cdot
\left( \matrix 1 & 0 \\ 0 & 0  \endmatrix
\right)\cdot
\left( \matrix 1 & r \\ 0 & 1  \endmatrix
\right)=
\left( \matrix 1 & r \\ \e(r) & \e(r)r  \endmatrix
\right)$$
for any $r\in R$ we deduce that
$$s\left( \matrix 1 & r \\ \e(r) & \e(r)r  \endmatrix
\right)=
\left( \matrix \id & [r] \\ 0 & \id  \endmatrix
\right)^*(p_1^*L(1)\ot p_2^*\eta(1))=
(p_1+[r]p_2)^*L(1)\ot p_2^*\eta(1).$$
Combining this with (\ref{s1}) and (\ref{s2})
one can easily compute that
\begin{equation}\label{s3}
s\left( \matrix 0 & r \\ \e(r) & 0  \endmatrix
\right)=(\phi\times [r])^*\PP\ot p_1^*\eta(-\e(r)r)\ot
p_2^*([r]^*L(1)\ot L(-\e(r)r)).
\end{equation}
Note that the formulas (\ref{s1}), (\ref{s2}), and (\ref{s3})
completely determine $s$. Now the identity
$$\left( \matrix 1 & \e(r) \\ 0 & 1  \endmatrix
\right)\cdot
\left( \matrix 0 & 0 \\ 0 & r_0  \endmatrix
\right)\cdot
\left( \matrix 1 & 0 \\ r & 1  \endmatrix
\right)=
\left( \matrix \e(r)r_0r & \e(r)r_0 \\ r_0r & r_0  \endmatrix
\right)$$
which holds for any $r\in R$, $r_0\in R^+$ implies
the equality
\begin{equation}\label{s4}
\left( \matrix {\id} & 0 \\ {[r]} & {\id}  \endmatrix
\right)^*(p_2^*L(r_0)\ot p_1^*\eta(r_0))=
s\left( \matrix {\e(r)r_0r} & {\e(r)r_0} \\ {r_0r} & {r_0}  \endmatrix
\right).
\end{equation}
Computing the right hand side using (\ref{s1})--(\ref{s3})
and restricting to $A\times 0$ we obtain the
identity
\begin{equation}\label{s5}
[r]^*L(r_0)=L(\e(r)r_0r)\ot\eta(-r_0r\e(r)r_0).
\end{equation}
On the other hand, setting $r=1$ and restricting (\ref{s4})
to $0\times A$ we get
\begin{equation}\label{s6}
L(r_0^2)=[r_0]^*L(1)\ot\eta(r_0).
\end{equation}
Setting $r=1$ in (\ref{s5}) we obtain the triviality
of $\eta(r_0^2)$. Then taking $r_0=1$ and $r\in R^+$
in (\ref{s5}) we obtain that
$$[r]^*L(1)=L(r^2)$$
for $r\in R^+$. Comparing this with (\ref{s6}) we deduce
the triviality of $\eta(r_0)$ for all $r_0\in R^+$.
Now (\ref{s5}) implies that $L(\cdot)$ gives a
$\Sigma_{R,\e}$-structure.
\ed

\begin{cor} A $\Sigma_{R,\e}$-structure for $\phi$
exists if and only if a splitting of the homomorphism
(\ref{homSL}) for $n=2$ exists.
\end{cor}

\begin{ex} Let $E=\C/\Z[i]$ be an elliptic curve with complex
multiplication by the ring of Gaussian numbers $R=\Z[i]$,
so that the corresponding Rosati involution $\e$ is just
the complex conjugation. In this situation there is no
$\Sigma_{R,\e}$-structure corresponding to the standard
polarization of $E$. Indeed, the corresponding line
bundle $L(1)$ should be of the form $\O(x)$ where $x$ is
a point of order 2 on $E$. Now the identity
$[1+i]^*L(1)=L(2)$   leads   to   a    contradiction    (see
\cite{thetaid} for details).
\end{ex}

\section{Abelian schemes with real multiplication}\label{real}
      
Let $F$ be a totally real number field,
$R$ be its ring of integers.
Let $A\ra S$ be an abelian  scheme  with   real
multiplication  by  $R$, i.e. a ring homomorphism
$R\ra\End_S(A):r\mapsto [r]_A$ is given. Then the dual
abelian scheme $\hat{A}$ also has a natural  real
multiplication by $R$.
Let  $J,  \hat{J}\sub  F$  be  fractional  ideals for $R$ (=
non-zero finitely generated $R$-submodules of $F$)
such that $J\hat{J}\sub R$.
                        
\begin{defi} An $(J, \hat{J})$-polarization on $A$ is a pair
of $R$-module homomorphisms
\begin{align*}
&\la_J:J\ra\Hom_{R}^{\sym}(A,\hat{A}),\\
&\la_{\hat{J}}:\hat{J}\ra\Hom_{R}^{\sym}(\hat{A},A)
\end{align*}
where   $\Hom_{R}^{\sym}(A,\hat{A})$  is  an  $R$-module  of
symmetric $R$-linear homomorphisms $f:A\ra\hat{A}$
(i.e. $\hat{f}=f$ and $f\circ [r]_A=[r]_{\hat{A}}\circ f$ for any  $r\in
R$), such  that $\la_{\hat{J}}(m)\circ\la_J(l)=[lm]_A$ and
$\la_J(l)\la_{\hat{J}}(m)=[lm]_{\hat{A}}$  for  any $l\in J$, $m\in \hat{J}$.
\end{defi}
   
\begin{rem} Usually one also imposes some positivity
condition on a polarization. In the case of $(J,\hat{J})$-polarizations
one can fix  an ordering on $J$: this means that for each
embedding $\si:F\ra\R$ an orientation of the line
$J\ot_{R,\si}\R$ is chosen. Then
one should require that if  an  element  $l\in  J$  is  totally
positive then the homomorphism $\la_J(l):A\ra\hat{A}$ is positive
(i.e. $\la_J(l)$ is a polarization in the classical sense).
\end{rem}
      
Note that the notion of $(J,\hat{J})$-polarization is equivalent to
that of $(Jx, \hat{J}x^{-1})$-polarization for any $x\in F^*$.
Also an $(J,\hat{J})$-polarization   of   $A$   is   the   same  as
an $(\hat{J},J)$-polarization of $\hat{A}$.
When  $\hat{J}=J^{-1}$  we  recover the notion  of
$J$-polarization in the sense of P.~Deligne and G.~Pappas \cite{DePa}
(except for the positivity condition).
Recall that they define
an $J$-polarization of  an  abelian  scheme  $A$  with  real
multiplication  by $R$ as an $R$-linear  homomorphism
$\la:J\ra\Hom_{R}^{\sym}(A,\hat{A})$ such that
the image of a totally positive element of $J$ under $\la$
is positive, and the induced morphism $A\ot_{R} J\ra\hat{A}$
is an isomorphism.
In  this  case we have also an isomorphism
$\hat{A}\ot_{R}  J^{-1}\wt{\ra}  A$   which   induces   an
$R$-linear                                    homomorphism
$\la':J^{-1}\ra\Hom_{R}^{\sym}(\hat{A},A)$,
hence we get an $(J,J^{-1})$-polarization in our sense.
Conversely, given an $(J,J^{-1})$-polarization as above
then the morphism
$\mu:A\ot_R J\ra\hat{A}$ induced by $\la_J$
and the morphism $\mu':\hat{A}\ra A\ot_R J$ induced by
$\la_{J^{-1}}$ are inverse to each other, so that $\la_J$
gives an $J$-polarization of $A$ (except for the positivity
condition).

For a pair $(J,\hat{J})$ as above we define the subgroup
$\Ga(J,\hat{J})\sub\SL_2(F)$ as follows:
$$\Ga(J,\hat{J})=\{
\left( \matrix a_{11} & a_{12}\\ a_{21} & a_{22}  \endmatrix
\right)\in\SL_2(F):  a_{11},  a_{22}\in  R,   a_{12}\in   J,
a_{21}\in \hat{J} \}.$$
Note that for $x\in F^*$ the homomorphisms
$\rho_{J,\hat{J}}$ and $\rho_{Jx,\hat{J}x^{-1}}$ are compatible with
the natural isomorphism $\Ga(J,\hat{J})\simeq \Ga(Jx,\hat{J}x^{-1})$
(induced by the conjugation by
$\left( \matrix x^{\frac{1}{2}} & 0\\ 0 & x^{-\frac{1}{2}} \endmatrix
\right)$).
In particular, in the case $R=\Z$ the group
$\Ga(J,\hat{J})$ is always isomorphic to the principal congruenz-subgroup
$\Ga_0(N)=\Ga(\Z,N\Z)\sub \SL_2(\Z)$ (for some $N>0$).
   
If $A$ is $(J,\hat{J})$-polarized then
using $\la_J$ and $\la_{\hat{J}}$ we can define a homomorphism
$$\rho_{J,\hat{J}}:\Ga(J,\hat{J})\ra\SL_2(A):
\left( \matrix a_{11} & a_{12}\\ a_{21} & a_{22}  \endmatrix
\right)\mapsto
\left( \matrix [a_{11}]_{\hat{A}} & \la_J(a_{12})\\
\la_{\hat{J}}(a_{21}) & [a_{22}]_A  \endmatrix
\right)$$

      
\begin{lem} For any non-zero element $r\in R$  (resp.  $l\in
J$, $m\in \hat{J}$)   the   corresponding   morphism  $[r]_A$ (resp.
$\la_J(l)$, $\la_{\hat{J}}(m)$) is an isogeny.
\end{lem}
      
\Pf . There exists a non-zero integer $N$ such that $r'=N/r\in R$, so
that $[r']_A\circ [r]_A=[N]_A$. Hence, $[r]_A$ is an isogeny on each
fiber, therefore, it is an isogeny. Similar argument works for
$\la_J(l)$ and $\la_{\hat{J}}(m)$.
\ed
      
Recall   that  we  denote  by  $\Ga(A)$  the  image  of  the
homomorphism $\widehat{\SL}_2(A)\ra\SL_2(A)$.
Let  assume  for  simplicity  that the image $\rho_{J,\hat{J}}$ is
contained in $\Ga(A)$. Otherwise, we can consider a
finite  flat  base  change  $S'\ra S$ such that the image of
$\Ga(J,\hat{J})$ is contained in $\Ga(A_{S'})$ where $A_{S'}$ is the induced
abelian scheme over $S'$. Indeed, recall that one has an exact
sequence
$$0\ra\hat{A}_2\ra\Pic^{\sym}(A)\ra\und{\Hom}^{\sym}(A,\hat{A})\ra0$$
of sheaves in fppf topology, hence the boundary homomorphism
$J\stackrel{\la_J}{\ra}\Hom^{\sym}(A,\hat{A})\ra H^1(S,\hat{A}_2)$ which can
be considered as a $J^*\ot\hat{A}_2$-torsor over $S$ where
$J^*=\Hom_\Z(J,\Z)$   (similarly,   for   $\la_{\hat{J}}$  we  get  an
$\hat{J}^*\ot A_2$-torsor). Now we can take $S'$ to be  the
corresponding $J^*\ot\hat{A}_2\times \hat{J}^*\ot A_2$-torsor
over $S$.
                                               
By the above lemma we can define an obstruction homomorphism
$\de_{J,\hat{J}}:\Ga(J,\hat{J})\ra\Br(S)$ as follows:
if $a_{21}$-entry of the matrix $h\in\Ga(J,\hat{J})$
is non-zero then the same entry of $\ov{g}=\rho_{J,\hat{J}}(h)$ is an
isogeny   and   we   can   put  $\de_{J,\hat{J}}(h)=e(G_g)$  where
$g\in\wt{\SL}_2(A)$   lies   above   $\ov{g}$.    Otherwise,
$\de_{J,\hat{J}}(h)=0$. As in proposition \ref{obshom} one can check that
$\de_{J,\hat{J}}$  is  a  homomorphism.
      
Let $I\sub R$ be a non-zero ideal.
Let us denote by $\ov{\Ga}_I(J,\hat{J})$ the group of matrices
$\left( \matrix \bar{a}_{11} & \bar{a}_{12}\\ \bar{a}_{21} &
\bar{a}_{22}  \endmatrix\right)$  where
$\bar{a}_{11},\bar{a}_{22}\in R/I$, $\bar{a}_{12}\in J/IJ$,
$\bar{a}_{21}\in  \hat{J}/I\hat{J}$,          such          that
$\bar{a}_{11}\bar{a}_{22}-\bar{a}_{12}\bar{a}_{21}=1$.
Here  we  use  the  natural
homomorphism of $R/I$-modules $J/IJ\ot \hat{J}/I\hat{J}\ra R/IR$ induced
by $J\ot \hat{J}\ra R$.
      
\begin{lem}\label{surj} Assume that $I\sub J\hat{J}$. Then
the natural reduction homomorphism
$\pi_I:\Ga(J,\hat{J})\ra\ov{\Ga}_I(J,\hat{J})$ is surjective.
\end{lem}
      
\Pf . Let
$\left( \matrix \bar{a}_{11} & \bar{a}_{12}\\ \bar{a}_{21} &
\bar{a}_{22}  \endmatrix\right)\in \ov{\Ga}_I(J,\hat{J})$   be  any  element.
Choose any non-zero liftings
$a_{12}\in J$, $a'_{21}\in \hat{J}$ and $a'_{22}\in  R$  of
$\bar{a}_{12}$, $\bar{a}_{21}$ and $\bar{a}_{22}$.
For  an  element $r\in R$ and a finite $R$-module $Q$ we say
that  $r$  is  relatively  prime  to  $Q$  if   $Q=rQ$,   or
equivalently,  $r\not\in\pg$  for  any  prime  ideal   $\pg$
associated with $Q$. By the  Chinese
remainder  theorem  we can lift $\bar{a}_{11}$ to an element
$a_{11}\in R$ which is relatively prime to $J/Ra_{12}$.
On the other hand, $a_{11}$ is relatively prime to $R/J\hat{J}$ since
$a_{11}a'_{22}\equiv 1\mod(I)$ and $I\sub J\hat{J}$.  Therefore,  $a_{11}$  is
relatively prime to $J/J\hat{J}a_{12}$ (since
$\supp(J/J\hat{J}a_{12})\sub\supp(R/J\hat{J})\cup\supp(J/Ra_{12})$).
It follows that $J=a_{11}J+J\hat{J}a_{12}=(Ra_{11}+\hat{J}a_{12})J$ which
implies the equality $R=Ra_{11}+\hat{J}a_{12}$. Thus,
we can write $1=ra_{11}+ma_{12}$ where $r\in R$, $m\in \hat{J}$.
Let $x=a_{11}a'_{22}-a_{12}a'_{21}-1\in I$.
Then replacing $a'_{22}$ and $a'_{21}$ by
$a_{22}=a'_{22}-xr$, $a_{21}=a'_{21}-xm$
we achieve $a_{11}a_{22}-a_{12}a_{21}=1$.
\ed
      
\begin{prop}\label{triv2}
The homomorphism $\de_{J,\hat{J}}:\Ga(J,\hat{J})\ra\Br(S)$ is trivial.
\end{prop}
      
\Pf .
Consider the reduction homomorphism
$\pi_I:\Ga(J,\hat{J})\ra\ov{\Ga}_I(J,\hat{J})$
where $I\sub 2J\hat{J}$.
Then by Lemma \ref{surj} $\pi_I$ is surjective.
On the other hand, the kernel of $\pi_I$ is contained in the subgroup
$\Ga(J,2\hat{J})\sub\Ga(J,\hat{J})$.
By Proposition \ref{triv1} the restriction of $\de_{J,\hat{J}}$ to
the subgroup $\Ga(J,2\hat{J})$ is trivial. Hence,
$\de_{J,\hat{J}}(\ker(\pi_I))=0$, so that
$\de_{J,\hat{J}}=\ov{\de}\circ\pi_I$ for some homomorphism
$\ov{\de}:\ov{\Ga}_I(J,\hat{J})\ra\Br(S)$. Moreover,
since by Lemma \ref{surj} the homomorphism
$\Ga(J,2\hat{J})\ra\ov{\Ga}_I(J,2\hat{J})$ is surjective,
it follows that $\ov{\de}$ is trivial on
matrices with $\bar{a}_{21}\in 2\hat{J}/I\hat{J}$.
In particular, $\ov{\de}$ vanishes on any diagonal matrix.
Let
$h=\left(\matrix \bar{a}_{11} & \bar{a}_{12}\\ \bar{a}_{21} &
\bar{a}_{22}  \endmatrix\right)$ be any element of
$\ov{\Ga}_I(J,\hat{J})$.   Then    $\bar{a}_{11}\bar{a}_{22}\equiv
1\mod  (J\hat{J})$,  hence  $\bar{a}_{11}\mod  (J\hat{J})$  is a unit in
$R/J\hat{J}$. Let $u\in (R/I)^*$ be any unit such that
$u\equiv a_{11}\mod (J\hat{J})$ (such a unit  always  exists  since
$R/I$ is an artinian ring). Then replacing $h$ by
$h\cdot\left(\matrix u & 0\\ 0 & u^{-1} \endmatrix\right)$
we reduce the problem of showing that $\ov{\de}(h)=0$
to the case when $\bar{a}_{11}\equiv 1\mod(J\hat{J})$.
      
Now we use
the result of L.~Vaserstein \cite{Vas} which asserts
that if $F\neq \Q$ then
the subgroup of $\Ga(J,\hat{J})$ consisting of matrices with
$a_{11}\equiv 1\mod(J\hat{J})$ is  generated  by  elementary
matrices, i.e. matrices of the form
$\left( \matrix 1 & l \\ 0 & 1 \endmatrix\right)$ and
$\left( \matrix 1 & 0 \\ m & 1 \endmatrix\right)$, where
$l\in L$, $m\in M$. In the case $F=\Q$ we may assume that
$J=\Z$, $\hat{J}=N\Z$ for some $N\in\Z$ and the corresponding
assertion for $\Ga(\Z,N\Z)=\Ga_0(N)$ is trivial.
Note that $\de_{J,\hat{J}}$ vanishes on
elementary  matrices  (the corresponding Heisenberg groups are
either Mumford groups of line bundles or trivial, so they
admit Shr\"odinger representations), hence it vanishes on any
matrix with $a_{11}\equiv 1\mod(J\hat{J})$ and we are done.
\ed
      
Let $\widehat{\Ga}(J,\hat{J})$ be the preimage of $\Ga(J,\hat{J})$ under the
homomorphism $\widehat{\SL}_2(A)\ra\SL_2(A)$.
      
\begin{thm} For every $(J,\hat{J})$-polarized  abelian  scheme  $A$
over $S$ with
$R$-multiplication such that the image of $\rho_{J,\hat{J}}$
is contained in $\Ga(A)$ there exists a faithful action of
a central extension of the group $\widehat{\Ga}(J,\hat{J})$
by $\Z\times\Pic(S)$ on $\D^b(A)$.
Without this assumption we always have compatible faithful projective
actions of $\Ga(2J,2\hat{J})$ on $\D^b(A)$ and of
$\widehat{\Ga}(J,\hat{J})$
on $A_{S'}$ for some finite flat base change $S'\ra S$.
\end{thm}

\begin{cor} Let $A$ be  an  abelian scheme over a normal
noetherian connected base $S$. Assume that the projection
$\widehat{\SL}_2(A)\ra\SL_2(A)$ is surjective and
$\End_S(A)\simeq R$ is a totally real field. Then  there  is  a
faithful
action of a central extension of $\widehat{\SL}_2(A)$ by
$\Z\times\Pic(S)$ on $\D^b(A)$.
\end{cor}
      
\Pf . By Prop. X 1.5 of \cite{Geer} the general fiber $A_K$ admits an
$R$-linear polarization $\la:A_K\ra\hat{A}_K$. Hence,
$J=\Hom_K(A_K,\hat{A}_K)$ and $\hat{J}=\Hom_K(\hat{A}_K,A_K)$ can
be  considered  as  fractional  ideals  for  $R$,  such that
$J\hat{J}\sub R$. By definition $\SL_2(A_K)=\Ga(J,\hat{J})$ and
by Lemma \ref{extend} we have an $(J,\hat{J})$-polarization
on $A$.
\ed

Now let us consider the case of abelian scheme $A$ with
$R$-linear   principal   polarization $\phi:A\wt{\ra}\hat{A}$.
In this case we have a natural inclusion
$$i_{\phi}:\Sp_{2n}(R)\ra\SL_2(A^n):
\left( \matrix M_{11} & M_{12}\\ M_{21} & M_{22}  \endmatrix
\right)\mapsto
\left( \matrix [M_{11}]_{\hat{A}} & \phi_{(n)}[M_{12}]_A\\
{[M_{21}]_{\hat{A}}\phi_{(n)}^{-1}} & [M_{22}]_A  \endmatrix
\right)$$
where $A^n$ is the relative $n$-th cartesian power of $A$
with the induced polarization $\phi_{(n)}$,
$M_{ij}\in\Mat_n(R)$, for every abelian scheme $A$ with
multiplication by $R$ we denote the natural map
$\Mat_n(R)\ra\End(A^n)$ by $M\mapsto [M]_A$.
  
Now we claim that if $R$ is unramified at $2$ then one can split
the extension $\widehat{\SL}_2(A^n)\ra\SL_2(A)$ over
$\Sp_{2n}(R)$ provided that a symmetric line bundle $L(1)$ on
$A$ is given such that $\phi_{L(1)}=\phi$.
According to Theorem \ref{Sig} it is sufficient to construct
a $\Sigma_{R,\id}$-structure for $\phi$.
Note that since $R$ is unramified at $2$
every element $r\in R$ can be represented in the form
$r=a^2+2b$ with $a,b\in R$. Now we define
$L(r)=[a]^*L(1)\ot (\phi,[b]_A)^*\PP$ where
$\PP$ is the Poincar\'e line bundle. It is easy to see
that $L(r)$ doesn't depend on a choice of $a$ and $b$,
and satisfies (\ref{CM1}) and (\ref{CM2}) with $\e=\id$.
The induced structure for $A^n$ and $\Mat_n(R)$ is given by
the homomorphism
$$\Mat^{\sym}_n(R)\ra\Pic^{\sym}(A^n):
B=(b_{ij})\mapsto L(B)=
\bigotimes_{i<j}(\phi p_i, [b_{ij}]_A p_j)^*\PP\ot
\bigotimes_i p_i^*L(b_{ii})$$
where  $\Mat^{\sym}_n(R)$  denotes  symmetric  matrices with
entries in $R$, $p_i:A^n\ra A$ is the
projection on the $i$-th factor. It is easy to see that
$\phi_{L(B)}=[B]_A$ and that
$[C]_A^*L(B)\simeq L(\sideset{^t}{}{C} B C)$ for any $C\in\Mat_n(R)$.
Now we can write the required splitting
$$\Sp_{2n}(R)\ra\widehat{\SL}_2(A^n):
M=\left( \matrix M_{11} & M_{12}\\ M_{21} & M_{22}  \endmatrix
\right)\mapsto (i_{\phi}(M),
(\phi^{-1}_{(n)})^*L(\sideset{^t}{_{11}}{M}M_{21}),
L(\sideset{^t}{_{22}}{M}M_{12})).$$

Using the above splitting  we  can  construct  a  projective
action of $\Sp_{2n}(R)$ on $D^b(A^n)$.
The vanishing of the obstacle
follows in this case from the  fact that $\Sp_{2n}(R)$ is
generated by elementary matrices established in \cite{BMS}.
      
\begin{thm}  Let  $A$  be  an  abelian  scheme   with   real
multiplication by $R$ over $S$, $L(1)$ be a
symmetric line bundle on $A$ rigidified along  the  zero
section  such  that  $\phi_{L(1)}:A\ra\hat{A}$ is an $R$-linear
isomorphism.
Assume  that  $R$  is  unramified at $2$. Then there is
a canonical faithful action of a central
extension of $\Sp_{2n}(R)$ by $\Z\times\Pic(S)$ on
$\D^b(A^n)$ where $A^n$ is the relative cartesian  power  of
$A$, $n\ge 1$. These actions are compatible via the natural
embeddings $\D^b(A^n)\ra\D^b(A^{n+1})$ and
$\Sp_{2n}(R)\ra\Sp_{2n+2}(R)$.
\end{thm}
      
\Pf . The same argument as in Proposition \ref{obshom} allows
to define an obstacle homomorphism $\de:\Sp_{2n}(R)\ra\Br(S)$
such that $\de(h)=0$ if and only if there exists a global
object in $\Isom_{H(X)}(\rho,\rho^h)$ where
$X=\hat{A}^n\times_S A^n$, $\rho$ is the representation of the
Heisenberg groupoid on $\D^b(A^n)$.
It is easy to check that $\de$ vanishes on elementary matrices,
hence it is zero.
\ed
             
\section{The central extension}
                                           
In this section we describe
explicitly the central extension of $\Sp_{2n}(\Z)$ by
$\Z\times\Pic(S)$ corresponding  to  the  projective action
defined  in the previous section.
  
We are going to use  a  presentation  of  $\Sp_{2n}(\Z)$  by
generators and relations borrowed from \cite{Cl}.
We always use the
standard symplectic basis $e_1,\ldots,e_n,f_1,\ldots f_n$ in
$\Z^{2n}$ such that $(e_i,f_j)=\de_{i,j}$).
First of all let us introduce
the relevant elementary matrices following the notation of \cite{Cl} 5.3.1.
Let $\S_{2n}$ be the set of pairs $(i,j)$ where
$1\le i,j\le 2n$ which
are not of the form $(2k-1,2k)$ or $(2k,2k-1)$. Then for for every
$(i,j)\in \S_{2n}$ we define an elementary
matrix $E_{ij}$ as follows:
$$E_{2k,2l}=
\left(\matrix 1 & 0 \\ \ga_{k,l} & 1\endmatrix\right),$$
$$E_{2k-1,2l-1}=
\left(\matrix 1 & -\ga_{k,l} \\ 0 & 1\endmatrix\right),$$
$$E_{2k-1,2l}=
\left(\matrix e_{kl} & 0 \\ 0 & e_{lk}^{-1}\endmatrix\right),$$
$$E_{2l,2k-1}=E_{2k-1,2l}$$
where $\ga_{kl}$ has zero $(\a,\b)$-entry unless $(\a,\b)=(k,l)$
or $(\a,\b)=(l,k)$, in the latter case $(\a,\b)$-entry is 1;
$e_{kl}$ for $k\neq l$ is the usual elementary matrix with units on
the diagonal and at $(k,l)$-entry and zeros elsewhere.
Now theorem 9.2.13 of \cite{Cl} asserts that for $n\ge 3$
the group $\Sp_{2n}(\Z)$ has a presentation consisting of
the generators $E_{ij}=E_{ji}$ (where $(i,j)\in\S_{2n}$)
subject to the relations
\begin{enumerate}
\item $[E_{ij},E_{kl}]=1$, if
$(i,k), (i,l), (j,k), (j,l)$ are in $\S_{2n}$
\item $[E_{ij},E_{kl}]=E_{il}$, if $(j,k)\not\in\S_{2n}$, $j$
is even, and $i$, $j$, $k$, and $l$ are distinct
\item $[E_{ij},E_{ki}]=E_{ii}^2$, if $(j,k)\not\in\S_{2n}$, $j$
is even, and $i$, $j$, and $k$ are distinct
\item $[E_{ii},E_{kl}]=E_{il}E_{ll}^{-1}$ if
$(i,k)\not\in\S_{2n}$, $i$ is even, and $i$, $k$, and $l$ are
distinct
\item $[E_{ii},E_{kl}]=E_{il}^{-1}E_{ll}^{-1}$ if
$(i,k)\not\in\S_{2n}$, $i$ is odd, and $i$, $k$, and $l$ are
distinct
\item $(E_{11}E_{22}E_{11})^4=1$.
\end{enumerate}
                                   
It is convenient to introduce also the symplectic matrix
$$\varphi=
\left(\matrix 0 & -1 \\ 1 & 0\endmatrix\right).$$
Then one has the following relations
\begin{equation}\label{conj1}
\varphi^{-1}E_{2k-1,2l-1}\varphi=E_{2k,2l}
\end{equation}
for all $1\le k,l\le n$,
\begin{equation}\label{conj2}
\varphi^{-1}E_{2k-1,2l}\varphi=E_{2l-1,2k}^{-1}
\end{equation}
for all $k\neq l$.
In particular, the group $\Sp_{2n}(\Z)$ is generated by
$\varphi$ and $E_{ij}$ with $i$  odd.  The  latter  set  of
generators is more convenient from the point of view of
our projective
representation on $\D^b(A^n)$ since the functors corresponding
to $E_{ij}$ with $i$ odd are very easy to describe
(see the proof of the theorem below).
                       
Let us denote by $\wt{\Sp}_{2n}(\Z)$ the group with generators
$E_{ij}=E_{ji}$ ($(i,j)\in\S_{2n}$) and one more generator
$\eps$ subject to relations (1)--(5) above, the
commutativity relation $[\eps,E_{ij}]=1$ for all
$(i,j)\in\S_{2n}$, and the modified relation (6)
$$(E_{11}E_{22}E_{11})^4=\eps.$$
  
Let $\pi:A\ra S$ be an abelian scheme
with a symmetric
line bundle $L$ (rigidified along the zero section)
which induces a principal polarization
$\phi:A\wt{\ra}\hat{A}$.
Let us also denote $\De(L)=2\pi_*L+e^*\om_{A/S}\in\Pic(S)$
where $\om_{A/S}$ is the relative canonical bundle.
It is known that
$4\cdot\De(L)=0$ (see e.~g. \cite{FC}, I, 5.1).
  
\begin{thm}\label{centrext}
Let $n\ge 3$.
The group $\wt{\Sp}_{2n}(\Z)$ is a central extension of
$\Sp_{2n}(\Z)$ by $\Z$.
The central extension of $\Sp_{2n}(\Z)$ by $\Z\times\Pic(S)$
corresponding to the projective action on $\D^b(A^n)$
is obtained from the $\wt{\Sp}_{2n}(\Z)$ by the
push-forward with respect to the homomorphism
$\Z\ra\Z\times\Pic(S):
1\mapsto (2g,2\De(L))$
\end{thm}
  
\Pf . Let us choose the intertwining functors corresponding
to the generators $\varphi$ and $E_{ij}$ (with $i$ odd)
in the following way. The functor corresponding to
$\varphi$ is the composition $\phi_{(n)}^{-1}\circ F_{A^n}$
where $F_{A^n}:\D^b(A^n)\ra\D^b(\hat{A}^n)$ is the Fourier-Mukai
transform. The functor corresponding to $E_{2k-1,2l-1}$
is simply tensor multiplication with the line bundle
$L(\ga_{k,l})$. Note that $L(\ga_{k,l})=(\phi p_k,p_l)^*\PP$
if $k\neq l$ while $L(\ga_{kk})=p_k^*L(1)$. At last the
functor corresponding to $E_{2k-1,2l}$ is $[e_{lk}]_A^*$.
We claim that these  functors  satisfy  all the relations
(1)-(5).
Let $P\sub\Sp_{2n}(\Z)$ be the subgroup of matrices of
the form
$\left(\matrix * & * \\ 0 & *\endmatrix\right)$.
Then there is an obvious action
of $P$ on $\D^b(A^n)$ such that the element
$\left(\matrix\sideset{^t}{}C^{-1} & 0 \\ 0 & C
\endmatrix\right)
\left(\matrix 1 & B \\ 0 & 1\endmatrix\right)$
where $C\in\GL_n(\Z)$, $B\in\Mat^{\sym}(n,\Z)$,
acts by the functor $[C^{-1}]_A^*\circ(\cdot \ot L(-B))$.
It is easy to see that our definition of the functors
corresponding to the generators $E_{ij}$ for $i$ odd is
compatible with this action. This means that all the
relations out of (1)--(5) which contain only these generators
are satisfied by the corresponding functors. Furthermore,
using the relations (\ref{conj1}), (\ref{conj2}) one can see
that all the relations out of (1)--(5) containing a generator
$E_{ij}$ with $i$ and $j$ of opposite  parity  in  the  left
hand side, are satisfied by our functors. Similarly, the
relation   (5)   follows   from   (4)  using  the  relations
(\ref{conj1}) and (\ref{conj2}).
It remains to check the relation (1) for $i$ and $j$ even, and
$k$ and $l$ odd, the relations (2), (3) for $i$ even and $l$
odd, and the relation (4) for $l$ odd. This can be done
directly applying the both sides of a relation to the object
$e_*\O_S\in\D^b(A^n)$. Thus, the relations (1)--(5) hold for
our functors. Now using that
$F_A^2\simeq [-1]^*(\cdot)\ot\om_{A/S}^{-1}[-g]$, where
$g=\dim A/S$, and that $F_A(L)\simeq L^{-1}\ot\pi^*\pi_*L$
one can easily compute that the functor corresponding
to $(E_{11}E_{22}E_{11})^4$ is
$(\cdot)\ot\De(L)^{\ot 2}[2g]$.
\ed
  
The $\Z$-part of the central extension of $\Sp_{2n}(\Z)$
acting on $\D^b(A^n)$ was computed in \cite{Orlov}.
Namely, for $g=1$ the corresponding class
in $H^2(\Sp_{2n}(\Z),\Z)$ is a half of the class of the
cocycle given by the Malsov index. 
On the other hand, it is easy to see that 
the class of the central extension
$\wt{\Sp}_{2n}(\Z)$ is a generator of $H^2(\Sp_{2n}(\Z),\Z)$
for sufficiently large $n$.
Indeed, it is known that $H^2(\Sp_{2n}(\Z),\Z)=\Z$ for large $n$
while $\Sp_{2n}(\Z)=[\Sp_{2n}(\Z),\Sp_{2n}(\Z)]$ for $n\ge 3$.
Moreover, the relations (2) and (4) easily imply that the element
$\eps$ belongs to $[\wt{\Sp}_{2n}(\Z),\wt{\Sp}_{2n}(\Z)]$, hence
$\wt{\Sp}_{2n}(\Z)=[\wt{\Sp}_{2n}(\Z),\wt{\Sp}_{2n}(\Z)]$.
It follows, that the central extension of $\Sp_{2n}(\Z)$ by
$\Z/p\Z$ obtained from $\wt{\Sp}_{2n}(\Z)$ is non-trivial for every
prime $p$, so our claim follows.
                                     
Let $\Ga_{1,2}\sub\Sp_{2n}(\Z)$ be the subgroup of matrices
$\left(\matrix M_{11} & M_{12} \\ M_{21} & M_{22}\endmatrix\right)$
such that $\sideset{^t}{_{11}}{M}M_{12}$ and
$\sideset{^t}{_{22}}{M}M_{21}$ have even diagonal entries.
  
Let $A$ be a principally polarized abelian scheme  over  $S$.
Then  we  have  a
canonical splitting of the projection
$\widehat{\SL}_2(A^n)\ra\SL_2(A^n)$
over  $\Ga_{1,2}$  which  is  constructed as in the previous
section using line bundles
$$L(B)=\bigotimes_{i<j}(\phi p_i, [b_{ij}]_A p_j)^*\PP\ot
\bigotimes_i (\phi p_i, [b_{ii}/2] p_i)^*\PP$$
associated with symmetric integer even-diagonal matrices
$B=(b_{ij})$ (note that this time we don't need any additional
data on $A$). It is known (see \cite{ThetaII} A.4)
that $\Ga_{1,2}$ is generated by
elements
$$\varphi,
\left(\matrix \sideset{^t}{^{-1}}{C} & 0 \\ 0 & C
\endmatrix\right),
\left(\matrix 1  &  B \\ 0  &  1 \endmatrix\right)$$
where $C\in\GL_n(\Z)$, $B$ is symmetric integer with even
diagonal. Obviously, this implies
vanishing of the obstruction for the projective action of
$\Ga_{1,2}$ on $\D^b(A^n)$ by intertwining operators,
hence this leads to a central extension of $\Ga_{1,2}$ by
$\Z\times\Pic(S)$.
  
\begin{prop} Let $A/S$  be  a  principally  polarized  abelian
scheme of dimension $g\ge 3$. Then the central extension of
$\Ga_{1,2}$ by $\Pic(S)$ acting on $\D^b(A^n)$ up to shifts is
trivial.
\end{prop}
  
\Pf . Considering a finite flat covering of $S$ corresponding
to a choice of a symmetric line bundle inducing
a principal polarization and using Theorem \ref{centrext}
one can see that the central extension in question is induced
by a central extension of $\Ga_{1,2}$ by
the torsion subgroup $\Pic(S)^{\tors}\sub\Pic(S)$.
Note that it is sufficient to prove our assertion in the case
when $A$ is the universal abelian scheme over the moduli
stack $\AA_g$ of principally polarized abelian schemes.
It remains to notice that $\Pic(\AA_g)^{\tors}=0$ since
$\Sp_{2g}(\Z)$ has no abelian quotients for $g\ge 3$
(this is deduced using the Kummer exact sequence --- see \cite{Mu}).
\ed
  
\begin{cor} The central extension of $\Sp_{2n}(\Z)$ by
$\Z/2\Z$ obtained by push-forward from $\wt{\Sp}_{2n}(\Z)$
has a splitting over $\Ga_{1,2}$.
\end{cor}

\end{document}